\documentclass[letter]{article}

\usepackage[latin1]{inputenc}
\usepackage[T1]{fontenc}
\usepackage[francais,english]{babel}
\usepackage{amsmath}
\usepackage{amssymb}

\usepackage{graphicx}
\usepackage{hyperref}
\usepackage{amscd}
\usepackage{mathrsfs}
\usepackage{comment}

\date{\today}

\author{
Karine \textsc{Beauchard}\footnote{
Centre de Math\'ematiques Laurent Schwartz, Ecole Polytechnique, 91128 Palaiseau Cedex, FRANCE.
email: Karine.Beauchard@math.polytechnique.fr}
\thanks{The author was partially supported by the ``Agence Nationale de la Recherche'' (ANR),
Projet Blanc EMAQS number ANR-2011-BS01-017-01}, 
Horst \textsc{Lange}\footnote{
Mathematisches Institut, Universit\"at K\"oln, K\"oln, NRW, Germany
},
Holger \textsc{Teismann}\footnote{
Department of Mathematics and Statistics, Acadia University, Wolfville, NS, Canada
email: hteisman@acadiau.ca}
\thanks{The author is supported by the Natural Sciences and Engineering Research Council of Canada (NSERC). This paper was finished while this author 
was visiting the Centre de Math\'ematiques Laurent Schwartz in the Fall of 2012. The author would 
like to thank the Centre for its hospitality  and financial support.}, 
}

\title{Local exact controllability of a 1D  Bose-Einstein condensate in a time-varying box}

\newtheorem{thm}{Theorem}
\newtheorem{prop}[thm]{Proposition}

\newtheorem{lem}[thm]{Lemma}
\newtheorem{rk}[thm]{Remark}
\newcommand{\R}{\mathbb{R}}
\newcommand{\C}{\mathbb{C}}

\newcommand{\Zn}{\Phi}
\newcommand{\Lop}{\mathcal{L}}
\newcommand{\Mop}{\mathcal{M}}

\newcommand{\st}{{2}}

\newcommand{\muz}{\mu}

\newcommand{\scat}{\kappa}

\newcommand{\scrap}[1]{}
\newcommand{\disp}{\displaystyle}
\renewcommand{\Re}{\textrm{Re}}
\renewcommand{\Im}{\textrm{Im}}

\begin{document}

\maketitle
\begin{abstract}
We consider a one--dimensional Bose-Einstein condensate in a infinite square-well (box) potential.
This is a nonlinear control system in which the state is the wave function
of the Bose Einstein condensate and the control is the length of the box.
We prove that local exact controllability around the ground state (associated
with a fixed length of the box) holds generically with respect to the chemical potential $\muz $;
i.e. up to an at most countable set of $\muz $--values. 
The proof relies on the linearization principle 
and the inverse mapping theorem, as well as ideas from analytic perturbation theory.
\end{abstract}

\textbf{Key words:} quantum systems, controllability of PDEs.


\section{Introduction} 

\subsection{Background and original problem}

Controlled manipulation of Bose Einstein condensates (BECs) is an important objective in quantum control theory. In this paper we consider a one-dimensional condensate 
 in a hard-wall trap (``condensate-in-a-box"), where the trap size (box length) is a time-dependent 
 function $L(\tau )$, which provides the control. 
 The model (see (\ref{GPE}) below) was first proposed by Band, Malomed, and Trippenbach \cite{BMT02} to study adiabaticity in
 a nonlinear quantum system.  More recently, the opposite regime, fast transitions (``shortcuts to adiabaticity"),
  has been investigated for BECs in box potentials \cite{theodorakis2009oscillations,del2012shortcuts}.  
 Condensates in a box trap have also been realized experimentally \cite{meyrath2005bose}, an achievement that 
 attracted considerable attention. 
Motivated by these developments, we study the controllability of the following  
system
\cite{BMT02}
\begin{equation} \label{GPE}
\left\lbrace \begin{array}{ll}
i \hbar \partial_\tau \Phi (\tau,z)
= - \frac{\hbar ^2}{2m} \partial_z^2 \Phi (\tau,z)
\mp  \scat |\Phi|^2 \Phi(\tau,z), & z \in (0,L(\tau)), \tau \in (0,\tau^*),\\
\Phi(\tau,0)=\Phi(\tau,L(\tau))=0, & \tau \in (0,\tau^*).
\end{array} \right. 
\end{equation}
Here
$\hbar$ is Planck's constant,
$m$ is the particle mass,
$\scat >0$  is a nonlinearity parameter derived from the scattering length and the particle number,
$\tau^*>0$ is a positive real number
and $L \in C^0([0,\tau^*],\mathbb{R}^*_+)$ is the length of the box.
The `$-$' sign in the right-hand side refers to the focusing case (attractive two-particle interaction), while the `$+$' sign refers to the defocusing one (repulsive interaction).
In this article, we will work with classical solutions (point-wise solutions) of the system (\ref{GPE}).\\

System (\ref{GPE}) is a \textit{nonlinear control system} in which
\begin{enumerate}
\item the \textit{state} is the wave function $\Phi(\tau,z)$, which is normalized
\begin{equation} \label{norm}
\int_{0}^{L(\tau)} |\Phi(\tau,z)|^2 dz = 1, \quad \forall \tau \in (0,\tau^*);
\end{equation}
\item the \textit{control} is the length $L$ of the box, with 
\begin{equation} \label{BC:L}
L(0)=L(\tau^*)=1.
\end{equation}
\end{enumerate}
This problem is a nonlinear variant of the control problem studied by K.~Beauchard in \cite{SchroLgVar}\footnote{The study of 
the controllability of \textit{nonlinear} Schr\"odinger equations was proposed by Zuazua  \cite{zuazua}.
}.

\subsection{Change of variables}

Following Band et al. \cite{BMT02} we introduce new variables,
\begin{equation} \label{CVAR}
t:=\frac{\hbar}{2m} \int _0^\tau \frac{ds}{L(s)^2}, \quad \quad
x:=\frac{z}{L(\tau)}, \quad \quad
\Phi(\tau,z)=\frac{\hbar}{\sqrt{2 \scat m} L(\tau)} \psi(t,x),
\end{equation}
to non-dimensionalize the problem and to transform it to the time-independent domain $(0,1)$. 
Then defining 
\begin{equation} \label{def:u}
u(t):=\frac{2m}{\hbar} \dot{L}(\tau) L(\tau)
\end{equation}
or, equivalently
\begin{equation} \label{def:L}
L(\tau)=\exp\left( \int_0^t u(s) ds \right),
\end{equation}
we obtain
\begin{equation} \label{syst}
\left\lbrace \begin{array}{ll}
i \partial_t \psi = - \partial_x^2 \psi
\mp |\psi|^2 \psi
+ i u(t) \partial_x[x\psi], & x \in (0,1), t \in (0,T),
\\
\psi(t,0)=\psi(t,1)=0, & t \in (0,T),
\end{array} \right.
\end{equation}
where
\begin{equation} \label{def:T}
T:=\int_0^{\tau^*} \frac{ds}{L(s)^2}.
\end{equation}
The system (\ref{syst}) is a control system in which
\begin{enumerate}
\item the state is $\psi$ with 
$$\|\psi(t)\|_{L^2(0,1)}=\|\psi(0)\|_{L^2(0,1)} e^{\frac{1}{2} \int_0^t u },$$
\item the control is the real valued function $u$.
\end{enumerate}
Note that the previous changes of variables impose  constraints on the control $u$.
Indeed, the requirement $L(0)=L(\tau^*)=1$, together with (\ref{def:L}) and (\ref{def:T}) impose 
$$\int_0^T u=0.$$

In this article, we will work with classical solutions of (\ref{syst}), that will provide classical solutions of (\ref{GPE}).

To ensure that the controllability of (\ref{syst}) gives the one of (\ref{GPE}),
we need the surjectivity of the map $L \mapsto u$, which is proved in the next proposition.

\begin{prop}  \label{ivp}
Let $T>0$, $u\in L^\infty(0,T;\R )$ extended by zero on $(-\infty,0) \cup (T,\infty)$ and such that $\int_0^T u(t)dt=0$.
The unique maximal solution of the Cauchy problem
\begin{equation} \label{CYpb:t(tau)}
\left\lbrace \begin{array}{l}
g'(\tau) =   \frac{\hbar}{2m} e^{- 2 \int _0^{g(\tau)} u(s) ds },
\\
g(0) = 0,
\end{array}\right.
\end{equation}
is defined for every $\tau \geqslant 0$, strictly increasing and satisfies 
\begin{equation} \label{limit_ivp}
\lim_{\tau \rightarrow + \infty} g(\tau)=+\infty.
\end{equation}
thus $\tau^*:=g^{-1}(T)$ is well defined.
Let  $L:[0,\infty) \to [0,\infty) $ be defined by 
\begin{equation} \label{Lt} 
L(\tau) :=  \exp \left( \int _0^{g(\tau)} u(s) ds \right). 
\end{equation}
Then, (\ref{BC:L}) and (\ref{def:u}) are satisfied.
\end{prop} 

\noindent \textbf{Proof of Proposition \ref{ivp}:}  The function
$F:\mathbb{R} \rightarrow \mathbb{R}$ defined by $F(x):=\frac{\hbar}{2m} e^{- 2 \int _0^{x} u(s) ds }$ 
is continuous, globally Lipschitz (because $u \in L^\infty$), and uniformly bounded. 
By Cauchy-Lipschitz (or Picard Lindelof) theorem, there exists a unique solution to (\ref{CYpb:t(tau)}),
defined for every $\tau \in [0,+\infty)$. It is strictly increasing on $[0,+\infty)$ because $g'>0$.
Now, we prove (\ref{limit_ivp}) by contradiction. We assume that $g(\tau) \leqslant T$ for every $\tau \in [0,+\infty)$. Then,
$$g'(\tau) \geqslant  \frac{\hbar}{2m} e^{- 2 \|u\|_\infty T }, \forall \tau \in (0,+\infty)$$
thus
$$T \geqslant g(\tau) \geqslant  \frac{ \tau \hbar}{2m} e^{- 2 \|u\|_\infty T }, \forall \tau \in (0,+\infty)$$
which is impossible. Therefore, there exists $\tau_1>0$ such that $g(\tau_1)>T$. Then, 
$g' \equiv \hbar/2m$ on $(\tau_1,\infty)$, which implies (\ref{limit_ivp}). 
The relation (\ref{BC:L}) is satisfied because $g(\tau^*)=T$ and $\int_0^T u=0$.
By integrating the first equality of (\ref{CYpb:t(tau)}) and using (\ref{Lt}), we get
$$g(\tau)= \frac{\hbar}{2m} \int_0^{\tau} \exp\left(- 2 \int _0^{g(s)} u \right) ds = \frac{\hbar}{2m} \int_0^\tau \frac{ds}{L(s)^2}.$$
Thanks to (\ref{Lt}) and (\ref{CYpb:t(tau)}), we have
$$\frac{2m}{\hbar} \dot{L}(\tau) L(\tau)
= \frac{2m}{\hbar} g'(\tau) u[g(\tau)] \exp\left( 2 \int_0^{g(\tau)} u \right)
= u[g(\tau)]$$
which proves (\ref{def:u}). \hfill $\Box$

\subsection{Main result}

We introduce the unitary $L^2((0,1),\mathbb{C})$ sphere $\mathcal{S}$,
the operator $A$ defined by
$$\begin{array}{ll}
D(A):=H^2 \cap H^1_0((0,1),\mathbb{C}), &
A \varphi := - \varphi'',
\end{array}$$
and the spaces
\begin{equation} \label{def:Hs}
H^s_{(0)}((0,1),\mathbb{C}):=D(A^{s/2}), \forall s >0.
\end{equation}
In particular,
$$H^3_{(0)}((0,1),\mathbb{C})=\{ \varphi \in H^3((0,1),\mathbb{C}) ; \varphi=\varphi''=0 \text{ at } x=0,1 \}.$$
We also introduce, for $T>0$, the space
$$\dot{H}^1_0((0,T),\mathbb{R}):=
\left\{ u \in H^1_0((0,T),\mathbb{R}) ; \int_0^T u(t) dt =0 \right\}.$$
For $\mu \in (\mp \pi^2,+\infty)$, we denote by $\phi_\mu$ the nonlinear ground state; i.e. the unique 
positive solution of the boundary value problem 
\begin{equation} \label{def:phimu}
\left\lbrace \begin{array}{ll}
\phi_\mu'' \pm \phi_\mu^3 = \pm \mu \phi_\mu, & x \in (0,1), \\
\phi_\mu(0)=\phi_\mu(1)=0.                & 
\end{array} \right.
\end{equation}
(See Section \ref{section:GS} for existence and properties of $\phi_\mu$).
Then the couple $(\psi_\mu(t,x):=\phi_\mu(x) e^{\pm i\mu t}, u \equiv 0)$ is a trajectory of (\ref{syst}).
The goal of this article is to prove the local exact controllability of system
(\ref{syst}) around this reference trajectory, for generic $\mu \in (\mp \pi^2,+\infty)$.

\begin{thm} \label{thm:Main}
Let $T>0$. There exists a countable set $J \subset (\mp \pi^2,+\infty)$ such that,
for every $\mu \in (\mp \pi^2,+\infty) \setminus J$, 
the system (\ref{syst}) is exactly controllable in time $T$, locally around the ground state; i.e., 
there exists $\delta=\delta(\mu,T)>0$ and a $C^1$--map
$$\Upsilon:  \mathcal{V}  \rightarrow  \dot{H}^1_0((0,T),\mathbb{R}),$$
where
$$\mathcal{V}:=\{ \psi_f \in H^3_{(0)}((0,1),\mathbb{C}) ;
\| \psi_f - \phi_\mu e^{\pm i \mu T } \|_{H^3_{(0)}} < \delta \text{ and } \|\psi_f\|_{L^2}=\|\phi_\mu\|_{L^2}  \},$$
such that, $\Upsilon(\phi_\mu e^{\pm i \mu T})=0$, and for every $\psi_f \in \mathcal{V}$, 
the solution of (\ref{syst}) associated with the control $u:=\Upsilon(\psi_f)$, and the initial condition
\begin{equation} \label{IC=GS}
\psi(0,x)=\phi_\mu(x), \quad x \in (0,1)
\end{equation}
is defined on $[0,T]$ and satisfies $\psi(T)=\psi_f$.
\end{thm}

\begin{rk} 
Note that by the time reversibility of the Schr\"odinger equation this result may be generalized to include 
arbitrary initial data $\psi(0,.)=\psi_0$, which are close enough to $\phi_\mu$ in $H^3_{(0)}((0,1),\mathbb{C})$.
\end{rk}

\subsection{Structure of this article}

The proof of Theorem \ref{thm:Main} relies on a linearization principle, which involves  
proving the controllability of the linear system that arises by linearizing (\ref{syst}) around the trajectory
$(\psi_\mu(t,x):=\phi_\mu(x) e^{\pm i\mu t}, u \equiv 0)$ and  applying the inverse 
mapping theorem. Accordingly, this article is organized as follows.

After stating the existence and uniqueness of the ground state (Section \ref{section:GS}), i.e. the positive 
solution $\phi_\mu$ of (\ref{def:phimu}), 
we study in Section \ref{sec:WP} the well--posedness of the Cauchy problem associated with (\ref{syst}).
The $C^1$--regularity of the end--point map is established in Section \ref{sec:C1}. 
Section \ref{sec:spectral} contains a detailed description of the spectral properties of the linearized system.
In Section \ref{sec:Cont_Lin}, we prove the controllability of the linearized system,
under appropriate assumptions \emph{\textbf{(A)}} and \emph{\textbf{(B)}},
which, in Section \ref{sec:Generic}, are shown to hold generically with respect to 
the chemical potential $\mu \in (\mp \pi^2,+\infty)$. 
Finally, in Section \ref{sec:proof}, we prove the main result of this article. The final section of the main part of the paper 
contains some concluding remarks and perspectives (Section \ref{sec:ccl}).

The main body of the article is followed by four appendices, containing proofs omitted in the main part of the paper 
to improve its readability. 
In Appendix A the proof of Proposition \ref{Prop:GS} (Section \ref{section:GS}) on the existence of ground states is provided.
The spectral properties of the linearization stated in Section \ref{sec:spectral} are established in Appendix B. 
Appendix C contains the proof of the analyticity of the spectrum.
Finally, Appendix D deals with trigonometric moment problems: a classical result, which is used in Section \ref{sec:Cont_Lin}, is recalled here.

\subsection{A brief review of infinite-dimensional bilinear control systems}

In this section we provide references to some of the pertinent literature. We do not, however, attempt 
a comprehensive review of the field, which is beyond the scope of the 
this paper\footnote{For (partial) reviews of  (linear and bilinear) control of Schr\"odinger equations, see for example 
\cite{zuazua,Teismann-et-al,JMC-book}. 
On the even broader subject of quantum control, several review papers and monographs are available; for a recent survey, see e.g. \cite{brif2012control} and the literature  (680 references!) therein.}. 

Early controllability results for Schr\"odinger equations with bilinear controls
were negative; see  \cite{HTC,mirrahimi2004,teismann2005} and in particular \cite{Turinici} obtained by Turinici as  
a corollary to a more general result by Ball, Marsden and Slemrod \cite{ball-marsden-slemrod}. 
Turinici's result was adapted to nonlinear Schr\"odinger equations 
by Illner, Lange and Teismann \cite{Teismann-et-al}. Because of these non--controllability properties,
bilinear Schr\"odinger equations were considered to be non--controllable for a long time.
However, some progress was eventually made and the question is now better understood.

Concerning exact controllability, 
local and almost global (between eigenstates) results for 1D models were obtained by Beauchard \cite{KB-JMPA,SchroLgVar}
and Coron and Beauchard  \cite{KB-JMC}, respectively. 
In \cite{KB-CL}, Beauchard and Laurent proposed important simplifications of the proofs
and dealt with nonlinear Schr\"odinger and wave equations with bilinear controls,
but in simpler configurations than in the present article.
In \cite{JMC-CRAS-Tmin}, Coron proved that a positive minimal time may be required for the local controllability of the 1D model.
This subject was studied further by Beauchard and Morancey  \cite{Tmin}, 
and by Beauchard for 1D wave equations \cite{ondes}. Exact controllability has also been studied in infinite time by 
Nersesyan and Nersisian  \cite{VNHN,VNHN2}.

As for approximate controllability, 
Mirrahimi and Beauchard \cite{KB-MM} proved global approximate controllability
in infinite time for a 1D model, and Mirrahimi obtained a similar result
for equations with continuous spectrum \cite{MM}.
Using adiabatic theory and intersection of eigenvalues in the space of controls,  
Boscain and Adami proved approximate controllability in finite time for 
particular models \cite{Boscain-Adami}. 
Approximate controllability, in finite time, for more general models, has been studied
by 3 teams, using different tools:
Boscain, Chambrion, Mason, Sigalotti \cite{Chambrion-et-al, Chambrion-et-al2, Chambrion-et-al3},  
used geometric control methods;
Nersesyan \cite{Nersesyan1,Nersesyan2}
used feedback controls and variational methods;
and Ervedoza and Puel \cite{Ervedoza-Puel}
considered a simplified model.

Moreover, optimal control problems have been investigated for Schr\"odinger 
equations with a nonlinearity of Hartree type
by Baudouin, Kavian, Puel  \cite{Baudouin,Baudouin-Kavian-Puel} 
and by Cances, Le Bris, Pilot  \cite{CLB-Cances-Pilot}.
Baudouin and Salomon studied an algorithm for the computation of optimal controls \cite{Baudouin-Salomon}.
The idea of  ``finite controllability of infinite-dimensional systems" was introduced by 
Bloch, Brockett, and Rangan \cite{bloch2010}. Finally, we mention that the somewhat related problem of bilinear wave equations
was considered by Khapalov \cite{KhapalovCOCV,KhapalovDCDS,KhapalovDCDIS}, who proves  global
approximate controllability to nonnegative equilibrium states.

\subsection{Notation}

If $X$ is a normed vector space, $x \in X$ and $R>0$,
$B_X(x,R):=\{ y \in X ; \|x-y\|<R\}$ denote the open ball with radius $R$ and
$\overline{B}_X(x,R):=\{ y \in X ; \|x-y\| \leqslant R \}$ denotes the closed ball with radius $R$.
Implicitly, functions take complex values, thus we write, for instance $H^1_0(0,1)$ instead of $H^1_0((0,1),\mathbb{C})$.
Otherwise we specify it and write, for example $L^2((0,T),\mathbb{R})$, $L^2((0,1),\mathbb{C}^2)$, etc.
We denote by $\langle.,.\rangle $ the (complex valued) scalar product in $L^2((0,1),\mathbb{C}^2)$
\begin{equation} \label{bracket}
\langle U,V\rangle  = \textstyle \left\langle  \left( \begin{array}{c} U^{(1)} \\ U^{(2)} 
 \end{array} \right) ,\left( \begin{array}{c} V^{(1)} \\ V^{(2)} 
 \end{array} \right) \right\rangle    
 = \disp \int _0^1 \left[ U^{(1)}(x)\overline{V^{(1)}(x)} + U^{(2)}(x)\overline{V^{(2)}(x)} \right] dx,
\end{equation}
and the (complex valued) scalar product in $L^2((0,1),\mathbb{C})$
$$\langle f , g \rangle = \int_0^1 f(x) \overline{g(x)} dx.$$
When the symbols '$\pm$' (resp. '$\mp$') are used, the upper symbol '$+$' (resp. '$-$') refers to the focusing case,
while the lower symbol '$-$' (resp. '$+$') refers to the defocussing one. 
This convention holds in all the article, with only one exception explained in Remark \ref{rk:exception}.

\section{Ground states}
\label{section:GS}

In this brief section we establish existence, uniqueness 
and some important properties of the 
positive solutions $\phi_\mu$ of (\ref{def:phimu}).
Proofs will be provided in Appendix A.

\begin{prop} \label{Prop:GS}
For every $\mu \in (\mp \pi^2,+\infty)$, there exists a unique positive solution $\phi_\mu \in H^3_{(0)}((0,1),\mathbb{R})$ of (\ref{def:phimu}).
Moreover,  the map $\mu \in (\mp \pi^2,+\infty) \rightarrow \phi_\mu \in L^2(0,1)$ is analytic and
\begin{equation} \label{convexite}
\langle \partial_\mu \phi_\mu , \phi_\mu  \rangle >0, \quad \forall \mu \in (\mp \pi^2,+\infty),
\end{equation} 
\begin{equation} \label{limit:phi_mu}
\|\phi_\mu\|_{L^\infty(0,1)} \underset{\mu \rightarrow \mp \pi^2}{\longrightarrow} 0.
\end{equation}
\end{prop}

\begin{rk}
$\phi_\mu$ is actually a smooth function (of $x$), but this property will not be used in this article.
\end{rk}

\begin{rk}
Property (\ref{convexite}) is known in the literature as  \textit{``convexity condition"} or \textit{``Vakhitov-Kolokolov condition"}  or \textit{``slope condition"}; it plays an important r\^ole in the stability of solitary-wave solutions. 
\end{rk}

\section{Well posedness}
\label{sec:WP}

\begin{prop} \label{Prop:WP:psi}
Let $\mu \in (\mp \pi^2,+\infty)$, $T>0$, $u \in \dot{H}^1_0((0,T),\mathbb{R})$.
There exists a unique (classical) solution
$\psi \in C^0([0,T],H^3 \cap H^1_0) \cap C^1([0,T],H^1_0)$ of (\ref{syst})(\ref{IC=GS}). 
Moreover $\psi(T) \in H^3_{(0)}(0,1)$ and
$$\| \psi(t)\|_{L^2(0,1)} = \|\phi_\mu\|_{L^2(0,1)} e^{\frac{1}{2} \int_0^t u(s) ds}, \quad \forall t \in [0,T].$$ 
\end{prop}

This statement will be proved by working on the auxiliary system
\begin{equation} \label{syst_xi}
\left\lbrace \begin{array}{ll}
i \partial_t \xi = - \partial_x^2 \xi
- w(t) |\xi|^2 \xi
+ v(t) x^2 \xi, & x \in (0,1), t \in (0,T),
\\
\xi(t,0)=\xi(t,1)=0, & t \in (0,T),
\end{array} \right.
\end{equation}
with
$$w(t):= \pm e^{\int_0^t u} \qquad
\text{  and  }  \qquad
v(t):=\frac{1}{4} (\dot{u}-u^2)(t),$$
that results from (\ref{syst}) and the relation
\begin{equation} \label{CVAR:psi->xi}
 \psi(t,x) = \xi(t,x) e^{\frac{i}{4} u(t) x^2+\frac{1}{2} \int_0^t u(s) ds}. 
\end{equation}
The following proposition ensures the local (in time) well posedness of the associated
Cauchy-problem when $v$ is small enough in $L^2$.

\begin{prop} \label{Prop:WP}
Let $R_0>0$ and $r>0$.
There exists $T=T(R_0,r)>0$ and $\delta>0$ such that,
for every $\xi_0 \in H^3_{(0)}(0,1)$ with $\|\xi_0\|_{H^3_{(0)}} < R_0$,
$w \in L^\infty((0,T),\mathbb{R})$ with $\|w\|_{L^\infty(0,T)}<r$,
and $v \in L^2((0,T),\mathbb{R})$ with $\|v\|_{L^2(0,T)}<\delta$,
there exists a unique (classical) solution $\xi \in C^0([0,T],H^3_{(0)}) \cap C^1([0,T],H^1_0)$
of the system (\ref{syst_xi}) with the initial condition 
\begin{equation} \label{IC:xi}
\xi(0,x)=\xi_0(x), \quad x \in (0,1).
\end{equation}
Moreover $\|\xi(t)\|_{L^2(0,1)}=\|\xi_0\|_{L^2}$, $\forall t \in [0,T]$.
\end{prop}

The following technical result, proved in \cite[Lemma 1]{KB-CL},
will be used in the proof of Proposition \ref{Prop:WP}.

\begin{lem} \label{Lem:G}
Let $T>0$ and $f \in L^2((0,T),H^3 \cap H^1_0(0,1))$. The function
$G:t \mapsto \int_0^t e^{iAs} f(s) ds$
belongs to $C^0([0,T],H^3_{(0)}(0,1))$ and
\begin{equation} \label{majo_Lem:G}
\|G\|_{L^\infty((0,T),H^3_{(0)})} \leqslant
c_1(T) \|f\|_{L^2((0,T),H^3 \cap H^1_0)}
\end{equation}
where the constants $c_1(T)$ are uniformly bounded for $T$ lying in bounded intervals.
\end{lem}

\noindent \textbf{Proof of Proposition \ref{Prop:WP}:} 
Let $c_1$ be the constant of Lemma \ref{Lem:G} associated to the value $T=1$.
We introduce constants $c_2, c_2', c_3>0$ such that, for every $z,\tilde{z} \in  H^3_{(0)}(0,1)$,
\begin{equation} \label{def:c23}
\begin{array}{c}
\| |z|^2 z \|_{H^3_{(0)}} \leqslant c_2 \|z\|_{H^3_{(0)}}^3, 
\\
\| |z|^2 z - |\tilde{z}|^2 \tilde{z} \|_{H^3_{(0)}} \leqslant c_2' \|z-\tilde{z}\|_{H^3_{(0)}} \max\{ \|z\|_{H^3_{(0)}}^2 ; \|\tilde{z}\|_{H^3_{(0)}}^2 \}, 
\\
\| x^2 z \|_{H^3} \leqslant c_3 \|z\|_{H^3_{(0)}}. 
\end{array}
\end{equation}
We define 
\begin{equation} \label{def:R_delta_T}
R:=3R_0, \quad \quad \delta:=\frac{1}{3c_1 c_3}, \quad \quad \text{ and } \quad \quad  T=T(R_0,r):= \min\left\{ 1 ; \frac{1}{3 c_2 r R^2} ; \frac{1}{2 r c_2' R^2} \right\}.
\end{equation}
Let $v \in L^2((0,T),\mathbb{R})$ with $\|v\|_{L^2}<\delta$ and $w \in L^\infty((0,T),\mathbb{R})$ with $\|w\|_{L^\infty}<r$.
We introduce the map
$$\begin{array}{|crcl}
F: & \overline{B}_{C^0([0,T],H^3_{(0)})}(0,R) & \rightarrow & C^0([0,T],H^3_{(0)}) \\
   &    \xi                                   & \rightarrow & F(\xi)
\end{array}$$
where
$$F(\xi)(t)=e^{-iAt}\xi_0 -i \int_0^t e^{-iA(t-s)} \Big[ - w(s) |\xi|^2 \xi(s) + v(s) x^2 \xi(s) \Big] ds, \forall t \in [0,T].$$
Lemma \ref{Lem:G} proves that $F$ takes values in $C^0([0,T],H^3_{(0)})$.

\emph{First step: We prove that $F$ maps $\overline{B}_{C^0([0,T],H^3_{(0)})}(0,R)$ into itself. }
Using (\ref{def:R_delta_T}), we get, for every $t \in [0,T]$,
$$\| e^{-iAt}\xi_0 \|_{H^3_{(0)}} = \| \xi_0 \|_{H^3_{(0)}} < R_0 = \frac{R}{3},$$
$$\left\| \int_0^t e^{-iA(t-s)}   w(s) |\xi|^2 \xi(s)  ds \right\|_{H^3_{(0)}}
\leqslant
 r \int_0^t \| |\xi|^2 \xi(s)\|_{H^3_{(0)}} ds
\leqslant 
 r T c_2 R^3
\leqslant
\frac{R}{3}.$$
By Lemma \ref{Lem:G} and (\ref{def:R_delta_T}) we also have, for every $t \in [0,T]$,
$$\left\|  \int_0^t e^{-iA(t-s)} v(s) x^2 \xi(s) ds \right\|_{H^3_{(0)}} 
\leqslant c_1 \|v\|_{L^2(0,t)} \|x^2 \xi\|_{L^\infty((0,t),H^3)} 
\leqslant c_1 c_3 \|v\|_{L^2(0,T)} R 
< \frac{R}{3}.$$

\emph{Second step: We prove that $F$ is a contraction of $\overline{B}_{C^0([0,T],H^3_{(0)})}(0,R)$.}
Working as in the first step, we get, for any $\xi_1, \xi_2 \in \overline{B}_{C^0([0,T],H^3_{(0)})}(0,R)$ the following estimates
$$\begin{array}{ll}
\left\| \int_0^t e^{-iA(t-s)}  w(s) [ |\xi_1|^2 \xi_1(s) - |\xi_2|^2 \xi_2(s) ]  ds \right\|_{H^3_{(0)}}
& \leqslant T r c_2' \| \xi_1-\xi_2 \|_{L^\infty(H^3_{(0)})} R^2 \\
& \leqslant \frac{\| \xi_1-\xi_2 \|_{L^\infty(H^3_{(0)})}}{2} , 
\end{array}$$
$$\begin{array}{ll}
\left\| \int_0^t e^{-iA(t-s)}  v(s) x^2(\xi_1-\xi_2)(s)  ds \right\|_{H^3_{(0)}}
& \leqslant c_1 c_3 \|v\|_{L^2(0,T)} \|\xi_1-\xi_2\|_{L^\infty(H^3_{(0)})} \\
& \leqslant \frac{\|\xi_1-\xi_2\|_{L^\infty(H^3_{(0)})}}{3},
\end{array}$$
where $L^\infty(H^3_{(0)})=L^\infty((0,T),H^3_{(0)})$.
\\

\emph{Third step: Conclusion.} By applying the Banach fixed point theorem to the map $F$,
we get a function $\xi \in \overline{B}_{C^0([0,T],H^3_{(0)})}(0,R)$ such that $F(\xi)=\xi$.
From this equality, we deduce that $\xi \in C^1([0,T],H^1_0)$ and that the first equality of (\ref{syst_xi})
holds in $H^1_0(0,1)$ for every $t \in [0,T]$. In particular, $\xi$ is a classical solution of the equation. $\Box$
\\

The following proposition ensures that maximal solutions of (\ref{syst_xi}) are global in time.

\begin{prop} \label{Prop:WP_global}
Let $T>0$, $\xi_0 \in H^3_{(0)}(0,1)$, $v \in L^2((0,T),\mathbb{R})$ and $w \in H^1((0,T),\mathbb{R})$.
There exists a unique (classical) solution $\xi \in C^0([0,T],H^3_{(0)}) \cap C^1([0,T],H^1_0)$ of the system (\ref{syst_xi})(\ref{IC:xi}). 
There exists $C=C(\|\xi_0\|_{H^3_{(0)}},\|v\|_{L^2(0,T)},\|w\|_{H^1(0,T)})>0$ such that
$$\| \xi \|_{L^\infty((0,T),H^3_{(0)})} \leqslant C.$$
Moreover $\|\xi(t)\|_{L^2}=\|\xi_0\|_{L^2}$, $\forall t \in [0,T]$.
\end{prop}

Then, Proposition \ref{Prop:WP:psi} follows from Proposition \ref{Prop:WP_global} and the change of variable (\ref{CVAR:psi->xi}).
\\

\noindent \textbf{Proof of Proposition \ref{Prop:WP_global}:} We extend $v$ by zero and $w$ by $w(T)$ on $(T,+\infty)$. 
Our goal is to prove the existence and uniqueness of a solution 
$\xi \in C^0([0,+\infty),H^3_{(0)}) \cap C^1([0,+\infty),H^1_0)$ of (\ref{syst_xi})(\ref{IC:xi}).
\\

\emph{First step: Maximal solution.} By Proposition \ref{Prop:WP}, there exists a unique local (in time) solution 
$\xi \in C^0([0,T_1],H^3_{(0)}) \cap C^1([0,T_1],H^1_0)$ of (\ref{syst_xi})(\ref{IC:xi}), for some time $T_1>0$. 
The uniqueness of  Proposition \ref{Prop:WP} and Zorn Lemma imply the existence of a unique maximal solution 
$\xi \in C^0([0,T^*),H^3_{(0)}) \cap C^1([0,T^*),H^1_0)$ of (\ref{syst_xi})(\ref{IC:xi}), for some time $T^* \in (0,+\infty]$.
Now, we prove by contradiction that $T^*=\infty$. We assume that $T^*<+\infty$.
\\

\emph{Second step: We prove that $\xi(t)$ is bounded in $H^1_0(0,1)$ uniformly with respect to $t \in [0,T^*)$.}
We recall that $\xi \in C^1([0,T^*),H^1_0)$, and the first equality of (\ref{syst_xi})
holds in $H^1_0(0,1)$ for every $t \in [0,T]$. Thus, the function 
$$J(t):=\int_0^1 \left( \frac{1}{2} |\partial_x \xi(t,x)|^2 - \frac{w(t)}{4} |\xi(t,x)|^4  \right) dx$$
satisfies
\begin{equation} \label{dJ/dt}
\frac{dJ}{dt}(t)=2v(t) \Im \left( \int_0^1 x \overline{\partial_x \xi(t,x)} \xi(t,x) dx \right)  - \frac{\dot{w}(t)}{4}  \|\xi(t)\|_{L^4}^4.
\end{equation}
We also recall the existence of a constant $\mathcal{C}>0$ such that (Galiardo-Nirenberg inequality \cite[p. 147]{Brezis})
$$\|f\|_{L^4(0,1)} \leqslant \mathcal{C} \|f\|_{L^2(0,1)}^{3/4} \|\partial_x f \|_{L^2(0,1)}^{1/4}, \quad \forall f \in H^1_0(0,1).$$
For every $t \in [0,T^*)$, we have
$$\begin{array}{ll} 
- \frac{w(t)}{4} \|\xi(t)\|_{L^4}^4 
   & \leqslant \frac{\mathcal{C}}{4} \|w\|_{L^\infty(0,T^*)} \|\xi(t)\|_{L^2}^{3} \|\partial_x \xi(t) \|_{L^2} 
\\ & \leqslant \frac{1}{4} \|\partial_x \xi(t) \|_{L^2}^2 + \frac{\mathcal{C}^2}{16} \|w\|_{L^\infty(0,T)}^2 \|\xi_0\|_{L^2}^6
\end{array}$$
thus
\begin{equation} \label{J/H1}
J(t) \geqslant \frac{1}{4} \|\partial_x \xi(t) \|_{L^2}^2 - \frac{\mathcal{C}^2}{16} \|w\|_{L^\infty(0,T^*)}^2 \|\xi_0\|_{L^2}^6 , 
\quad \forall t \in [0,T^*).
\end{equation}
We deduce that
\begin{equation} \label{dJ/dt_1}
\begin{array}{ll}
2v(t) \Im \left( \int_0^1 x \overline{\partial_x \xi(t,x)} \xi(t,x) dx \right) 
& \leqslant 2 |v(t)| \|\partial_x \xi(t)\|_{L^2} \\
& \leqslant 4v(t)^2 + \frac{1}{4}\|\partial_x \xi(t)\|_{L^2}^2 \\
& \leqslant 4v(t)^2 + J(t) + \frac{\mathcal{C}^2}{16} \|w\|_{L^\infty(0,T^*)}^2 \|\xi_0\|_{L^2}^6
\end{array}
\end{equation}
and
\begin{equation} \label{dJ/dt_2}
\begin{array}{ll}
- \frac{\dot{w}(t)}{4}  \|\xi(t)\|_{L^4}^4
& \leqslant \frac{\mathcal{C}}{4} |\dot{w}(t)| \|\partial_x \xi(t)\|_{L^2} \\
& \leqslant \frac{\mathcal{C}^2}{16} \dot{w}(t)^2 + \frac{1}{4} \|\partial_x \xi(t)\|_{L^2}^2 \\
& \leqslant \frac{\mathcal{C}^2}{16} \dot{w}(t)^2 + J(t) + \frac{\mathcal{C}^2}{16} \|w\|_{L^\infty(0,T^*)}^2 \|\xi_0\|_{L^2}^6.
\end{array}
\end{equation}
From (\ref{dJ/dt}), (\ref{dJ/dt_1}), (\ref{dJ/dt_2}) and Gronwall lemma, we get
$$J(t) \leqslant \left( J(0)+\int_0^t \Big( 4v(s)^2 + \frac{\mathcal{C}^2}{16} [ \dot{w}(s)^2 + 2 \|w\|_{L^\infty(0,T^*)}^2 \|\xi_0\|_{L^2}^6 ]  \Big) ds    \right) e^{2t}, \forall t \in [0,T^*).$$
Thus, $J$ is bounded uniformly with respect to $t \in [0,T^*)$, and so is $\|\xi(t)\|_{H^1}$ (see (\ref{J/H1})).
\\

\emph{Third step: We prove that $\xi(t)$ is bounded in $H^3_{(0)}(0,1)$ uniformly with respect to $t \in [0,T^*)$.}
First, we recall the existence of a constant $\mathcal{C}$ such that 
$$\| |\xi|^2 \xi \|_{H^3_{(0)}} \leqslant \mathcal{C} \|\xi\|_{H^3_{(0)}} \|\xi\|_{H^1_0}^2, \quad \forall \xi \in H^3_{(0)}(0,1).$$
This follows from the explicit expression of $\partial_x^3[|\xi|^2\xi]$ and the Galiardo-Nirenberg inequality.
From the relation $\xi=F(\xi)$ in $C^0([0,T],H^3_{(0)})$ and Lemma \ref{Lem:G}, we get, for every $t \in [0,T^*)$,
$$\|\xi(t)\|_{H^3_{(0)}} \leqslant \|\xi_0\|_{H^3_{(0)}}+
\int_0^t |w(s)| \mathcal{C} \|\xi\|_{L^\infty(H^1)}^2 \|\xi(s)\|_{H^3_{(0)}} ds +
c_1(T^*) \left( \int_0^t |v(s)|^2 c_3^2 \|\xi(s)\|_{H^3_{(0)}}^2 ds \right)^{1/2}$$
(see (\ref{def:c23}) for the definition of $c_3$).
Using Cauchy-Schwarz inequality, we get
$$\begin{array}{ll}
\|\xi(t)\|_{H^3_{(0)}}^2  \leqslant 
& 3 \|\xi_0\|_{H^3_{(0)}}^2 + 3 t \int_0^t |w(s)|^2 \mathcal{C}^2 \|\xi\|_{L^\infty(H^1)}^4 \|\xi(s)\|_{H^3_{(0)}}^2 ds 
\\  & + 3 c_1(T^*)^2  \int_0^t |v(s)|^2 c_3^2 \|\xi(s)\|_{H^3_{(0)}}^2 ds.
\end{array}$$
Then Gronwall lemma proves that $\xi(t)$ is bounded in $H^3_{(0)}$ uniformly with respect to $t \in [0,T^*]$.
\\

\emph{Fourth step: Conclusion.} From the relation $\xi(t)=F(\xi)(t)$ and the third step, 
$\xi(t)$ satisfies the Cauchy-criterion in $H^3_{(0)}(0,1)$ when $[t \rightarrow T^*]$. Thus the maximal solution may be extended
after $T^*$, which is a contradiction. Therefore $T^*=+\infty$. $\Box$

\section{$C^1$-regularity of the end-point map}
\label{sec:C1}

By Proposition \ref{Prop:WP:psi}, we can consider, for any $T>0$ and $\mu \in (\mp \pi^2,+\infty)$ the end point map
$$\begin{array}{|cccc}
\Theta_{T,\mu}: & \dot{H}^1_0((0,T),\mathbb{R}) & \rightarrow & H^3_{(0)}(0,1) \cap \mathcal{S}_{\|\phi_\mu\|_{L^2}} \\
                &      u                        & \mapsto     & \psi(T) 
\end{array}$$
where $\psi$ is the solution of (\ref{syst})(\ref{IC=GS}) and
$\mathcal{S}_{\|\phi_\mu\|_{L^2}}$ is the $L^2((0,1),\mathbb{C})$-sphere with radius $\|\phi_\mu\|_{L^2}$.
The goal of this section is the proof of the $C^1$-regularity of $\Theta_{T,\mu}$.

\begin{prop} \label{Prop:C1}
Let $\mu \in (\mp \pi^2,+\infty)$ and $T>0$. 
The map $\Theta_{T,\mu}$ is $C^1$, moreover, for every $u,U \in  \dot{H}^1_0((0,T),\mathbb{R})$, we have $d\Theta_{T,\mu}(u).U=\Psi(T)$
where $\Psi$ is the solution of the linearized system
\begin{equation} \label{syst_lin}
\left\lbrace \begin{array}{ll}
i \partial_t \Psi = - \partial_x^2 \Psi
\mp [ 2 |\psi|^2 \Psi + \psi^2 \overline{\Psi}]
+ i U(t) \partial_x[x\psi], & x \in (0,1), t \in (0,T),
\\
\Psi(t,0)=\Psi(t,1)=0, & t \in (0,T).
\\
\Psi(0,x)=0,           & x \in (0,1).
\end{array} \right.
\end{equation}
and $\psi$ is the solution of (\ref{syst})(\ref{IC=GS}).
\end{prop}

This proposition will be proved by working first on the auxiliary system (\ref{syst_xi}).

\subsection{For the auxiliary system (\ref{syst_xi})}

For  $\mu \in (\mp \pi^2,+\infty)$, we introduce the end-point map of the auxiliary system
$$\begin{array}{|cccc}
\Omega_{T,\mu}: & L^2 \times H^1 ((0,T),\mathbb{R})& \rightarrow & H^3_{(0)}(0,1) \cap \mathcal{S}_{\|\phi_\mu\|_{L^2}} \\
                &      (v,w)                       & \mapsto     & \xi(T)
\end{array}$$
where $\xi$ is the solution (\ref{syst_xi}) with the initial condition 
\begin{equation} \label{IC:xi_GS}
\xi(0,x)=\phi_\mu(x), \quad x \in (0,1).
\end{equation}

\begin{prop} \label{Prop:C1_xi}
Let $T>0$. The map $\Omega_{T,\mu}$ is $C^1$, moreover, for every $(v,w) \in  L^2 \times H^1((0,T),\mathbb{R})$, 
we have $d\Omega_{T,\mu}(v,w).(V,W)=\zeta(T)$ where $\zeta$ is the solution of the linearized system
\begin{equation} \label{syst_lin_xi}
\left\lbrace \begin{array}{ll}
i \partial_t \zeta = - \partial_x^2 \zeta
- w(t) [ 2 |\xi|^2 \zeta + \xi^2 \overline{\zeta}]
+ v(t) x^2 \zeta
- W(t) |\xi|^2 \xi
+ V(t) x^2 \xi, & x \in (0,1), t \in (0,T),
\\
\zeta(t,0)=\zeta(t,1)=0,        & t \in (0,T).
\\
\zeta(0,x)=0,                   & x \in (0,1).
\end{array} \right.
\end{equation}
and $\xi$ is the solution of (\ref{syst_xi})(\ref{IC:xi_GS}).
\end{prop}

\noindent \textbf{Proof of Proposition \ref{Prop:C1_xi}:}

\emph{First step: Well posedness of (\ref{syst_lin_xi}).} 
Let $(v,w), (V,W) \in  L^2 \times H^1 ((0,T),\mathbb{R})$ and $\xi$ be the solution of (\ref{syst_xi})(\ref{IC:xi_GS}).
The well posedness of (\ref{syst_lin_xi}) may be proved with a fixed point argument in $C^0([0,T_1],H^3_{(0)})$
under a smallness assumption on $\|w\|_{L^1(0,T_1)}$ and $\|v\|_{L^2(0,T_1)}$ for the map to be contracting.
Then, iterating this argument on a finite number of intervals $[0,T_1]$, $[T_1,T_2]$,... we get 
the well posedness of (\ref{syst_lin_xi}) on the whole interval $[0,T]$.
\\

\emph{Second step: Local Lipschitz regularity of $\Omega_{T,\mu}$.}
Let $(v,w) \in  L^2 \times H^1 ((0,T),\mathbb{R})$ and $\xi$ be the solution of (\ref{syst_xi})(\ref{IC:xi_GS}).
Let $(V,W) \in  L^2 \times H^1 ((0,T),\mathbb{R})$ with $\|(V,W)\|_{L^2 \times H^1(0,T)} \leqslant 1$ and $\widetilde{\xi}$ be the solution of
$$\left\lbrace \begin{array}{ll}
i \partial_t \widetilde{\xi} = - \partial_x^2 \widetilde{\xi}
- (w+W)(t) |\widetilde{\xi}|^2 \widetilde{\xi}
+ (v+V)(t) x^2 \widetilde{\xi}, & x \in (0,1), t \in (0,T),
\\
\widetilde{\xi}(t,0)=\widetilde{\xi}(t,1)=0, & t \in (0,T),
\\
\widetilde{\xi}(0,x)=\phi_\mu(x),          & x \in (0,1).
\end{array} \right.$$
We claim that there exists a constant $C_1=C_1(\|v\|_{L^2},\|w\|_{H^1})>0$ (independent of $(V,W)$) such that
\begin{equation} \label{Omega_lipschitz}
\| \widetilde{\xi} - \xi \|_{L^\infty(H^3_0)} \leqslant C_1 \|(V,W)\|_{L^2 \times H^1}.
\end{equation}
By Proposition \ref{Prop:WP_global}, there exists $R=R(\|v\|_{L^2},\|w\|_{H^1})>0$ (independent of $V$ and $W$) such that
\begin{equation} \label{bound:R}
\|\xi\|_{L^\infty(H^3_{(0)})}, \| \widetilde{\xi} \|_{L^\infty(H^3_{(0)})}  \leqslant R.
\end{equation}
Thus, there exists $C_2=C_2(R)>0$ such that
$$\| |\widetilde{\xi}|^2 \widetilde{\xi} - |\xi|^2 \xi \|_{L^\infty(H^3_{(0)})} \leqslant C_2 \| \widetilde{\xi} -  \xi \|_{L^\infty(H^3_{(0)})}.$$
From the relation
$$\begin{array}{ll}
(\widetilde{\xi}-\xi)(t)= -i \int_0^t e^{-iAs} \left[
- w [ |\widetilde{\xi}|^2 \widetilde{\xi} - |\xi|^2 \xi ]  
- W |\widetilde{\xi}|^2 \widetilde{\xi} 
+ v x^2 (\widetilde{\xi} - \xi)  + V x^2 \widetilde{\xi}
\right](s) ds,
\end{array}$$
Lemma \ref{Lem:G} and (\ref{def:c23}), we get
$$\begin{array}{ll}
\|(\widetilde{\xi}-\xi)(t)\|_{H^3_{(0)}} \leqslant 
   & \int_0^t \left( |w(s)| C_2 \|(\widetilde{\xi}-\xi)(s)\|_{H^3_{(0)}} + |W(s)| c_2 R^3 \right) ds
\\ & + c_1(T) \left( \int_0^t \Big[ |v(s)|^2 c_3^2 \|(\widetilde{\xi}-\xi)(s)\|_{H^3_{(0)}}^2 + |V(s)|^2 c_3^2 R^2  \Big] ds  \right)^{1/2}.
\end{array}$$
Thus,
$$\begin{array}{ll}
\|(\widetilde{\xi}-\xi)(t)\|_{H^3_{(0)}}^2 \leqslant 
   & 4t \int_0^t \left( |w(s)|^2 C_2^2 \|(\widetilde{\xi}-\xi)(s)\|_{H^3_{(0)}}^2 + |W(s)|^2 c_2^2 R^6 \right) ds
\\ & + 2 c_1(T)^2 \int_0^t \Big[ |v(s)|^2 c_3^2 \|(\widetilde{\xi}-\xi)(s)\|_{H^3_{(0)}}^2 + |V(s)|^2 c_3^2 R^2  \Big] ds  
\end{array}$$
and we get (\ref{Omega_lipschitz}) thanks to Gronwall lemma.
\\
\\

\emph{Third step: Existence of a constant $C=C(\|v\|_{L^2},\|w\|_{H^1})>0$ such that}
$$\| \widetilde{\xi} - \xi - \zeta \|_{L^\infty(H^3_{(0)})} 
\leqslant C \|(V,W)\|_{L^2 \times H^1 }^2, \text{ when } \|(V,W)\|_{L^2 \times H^1} \leqslant 1.$$
Thanks to (\ref{bound:R}), there exists a constant $C_3=C_3(R)>0$ such that
$$\| |\widetilde{\xi}|^2 \widetilde{\xi} - |\xi|^2 \xi - 2|\xi|^2 (\widetilde{\xi}-\xi) - \xi^2 \overline{(\widetilde{\xi}-\xi)} \|_{L^\infty(H^3_{(0)})}
\leqslant C_3 \|  \widetilde{\xi} -  \xi \|_{L^\infty(H^3_{(0)})}^2.$$
Let $\Delta:=\widetilde{\xi}-\xi-\zeta$. From the relation
$$\begin{array}{ll}
\Delta(t)= -i \int_0^t e^{-iAs} \left[ \right. 
& 
- w(s) [ |\widetilde{\xi}|^2 \widetilde{\xi}(s)- |\xi|^2 \xi(s) - 2|\xi|^2 (\widetilde{\xi}-\xi)(s) - \xi^2 \overline{(\widetilde{\xi}-\xi)}(s) ]
\\ &
- w(s)[ 2|\xi|^2 (\widetilde{\xi}-\xi-\zeta) + \xi^2 \overline{(\widetilde{\xi}-\xi-\zeta)}
\\ &
- W(s) [ |\widetilde{\xi}|^2 \widetilde{\xi}(s)- |\xi|^2 \xi(s) ]
\\ &
+ v(s) x^2 \Delta(s) + V(s) x^2 (\widetilde{\xi}-\xi)(s) \left. \right] ds
\end{array}$$
we deduce that
$$\begin{array}{ll}
\|\Delta(t)\|_{H^3_0} \leqslant
& 
\int_0^t |w(s)| \Big( C_3 C_1^2 \|(V,W)\|_{L^2 \times H^1}^2 + 3 R^2 \|\Delta(s)\|_{H^3_{(0)}} \Big)  ds 
\\
& + \int_0^t   |W(s)| C_2 C_1 \|(V,W)\|_{L^2 \times H^1}  ds 
\\ &
+ \left(\int_0^t \Big[ |v(s)|^2 c_3^2 \|\Delta(s)\|_{H^3_{(0)}}^2 + |V(s)|^2 c_3^2 C_1^2 \|(V,W)\|_{L^2 \times H^1}^2 \Big] ds  \right)^{1/2}.
\end{array}$$
We conclude the proof by taking the square of this inequality and applying Gronwall lemma. \hfill $\Box$

\subsection{For the system (\ref{syst})}

We now prove Proposition \ref{Prop:C1}. 
First, we recall that, for every $u\in \dot{H}^1_0((0,T),\mathbb{R})$,  
$\Theta_{T,\mu}(u)=\Omega_{T,\mu}(v,w)$, where $w(t):= \pm e^{\int_0^t u}$ and $v(t):=\frac{(\dot{u}-u^2)(t)}{4} $.
Thus $\Theta_{T,\mu}$ is $C^1$ and
$$d\Theta_{T,\mu}(u).U = d \Omega_{T,\mu} (v,w).(V,W)
\quad \text{ where } \quad
V:=\frac{\dot{U}-2uU}{4}
\text{ and }
W:=\pm \left( \int_0^t U \right)  e^{\int_0^t u}.$$
This gives the conclusion because
$$\Psi(t,x)=
\left[ \zeta(t,x) + \left( \frac{i}{4} U(t) x^2 + \frac{1}{2} \int_0^t U(s) ds \right) \xi(t,x) \right] 
e^{\frac{i}{4} u(t) x^2+\frac{1}{2} \int_0^t u(s) ds}.$$

\section{Spectral analysis and consequences}
\label{sec:spectral}

In this section, we are interested in the linearized system around
the nonlinear trajectory $(\psi_\mu(t,x)=\phi_\mu(x) e^{\pm i\mu t},u=0)$ where
$\phi_\mu$ is defined by (\ref{def:phimu}), for $\mu \in(\mp \pi^2,+\infty)$,
$$\left\lbrace \begin{array}{ll}
i \partial_t \Psi = - \partial_x^2 \Psi
\mp  [ 2 |\psi_\mu|^2 \Psi + \psi_\mu^2 \overline{\Psi}]
+ i U(t) \partial_x[x\psi_\mu], & x \in (0,1), t \in (0,T),
\\
\Psi(t,0)=\Psi(t,1)=0, & t \in (0,T),
\\
\Psi(0,x)=0,           & x \in (0,1).
\end{array} \right.$$
As usual, the time dependence of the second term in the right hand side is eliminated by the  
transformation  
$$\Psi(t,x)=\widetilde{\Psi}(t,x) e^{\pm i \mu t}$$
which leads to
\begin{equation} \label{syst_lin_static}
\left\lbrace \begin{array}{ll}
i \partial_t \widetilde{\Psi} = - \partial_x^2 \widetilde{\Psi} 
\pm \mu \widetilde{\Psi}
\mp  [ 2 \phi_\mu^2 \widetilde{\Psi} + \phi_\mu^2 \overline{\widetilde{\Psi}}]
+ i U(t) (x\phi_\mu)', & x \in (0,1), t \in (0,T),
\\
\widetilde{\Psi}(t,0)=\widetilde{\Psi}(t,1)=0, & t \in (0,T),
\\
\widetilde{\Psi}(0,x)=0,           & x \in (0,1).
\end{array} \right.
\end{equation}
In this section, we will work with the real $(2\times 2)$-system arising from 
this equation, by decomposition in real and imaginary parts. Consider the matrix operator 
\begin{equation} \label{L}
\Lop_\mu:= 
\left( \begin{array}{cc} 
0        & L_\mu^- \\
-L_\mu^+ & 0 
\end{array} \right) 
\quad \text{ where } \quad
\left\lbrace \begin{array}{l}
L_\mu^- :=  - \Delta \pm \mu \mp  \phi_\mu^2, \\
L_\mu^+ :=  - \Delta \pm \mu \mp 3\phi_\mu^2.
\end{array}\right.
\end{equation} 
The previous equation takes the form
\begin{equation} \label{Z} 
\left\lbrace \begin{array}{l}
\partial_t Z = \Lop_\mu Z + U(t)   (x\phi_\mu)'
\left( \begin{array}{r} 
1 \\ 0
\end{array} \right), 
\\
Z(t,0)=Z(t,1)=0,
\\
Z(0,x)=0,
\end{array}\right.
\end{equation} 
where 
 \[ Z(t,x) := {\Re[\widetilde{\Psi}(t,x)] \choose \Im[\widetilde{\Psi}(t,x)]} .\] 
For convenience, we also define
\begin{equation} \label{L}
\Lop_{\mp \pi^2}:= 
\left( \begin{array}{cc} 
0 & -\Delta -\pi^2 \\
\Delta +\pi^2 & 0 
\end{array} \right), \qquad
\phi_{\mp \pi^2}:=0.
\end{equation}

The goal of this section is to establish the spectral properties of the operators $\Lop_\mu$
needed in the proof of the controllability of the linear system (\ref{syst_lin_static}) in Section \ref{sec:Cont_Lin}.

\subsection{Auxiliary operators}

It will be convenient to employ a similarity transformation (see \cite[(12.15)]{RSS05}). Let
\begin{equation} \label{def:J}
J:= \left( \begin{array}{cr} 1 & i\\
 1 & -i \end{array} \right).  
\end{equation}
Then, for any $\mu \in [\mp \pi^2,+\infty)$, we have 
\begin{equation} \label{DM}  
i\Lop_\mu  =  J^{-1}   \Mop_\mu J
\text{ where }
\Mop_\mu  :=  
\left( \begin{array}{rr} 
-\Delta   &  0     \\
 0        & \Delta     
\end{array} \right) 
+  
\left( \begin{array}{cc} 
\pm \mu \mp 2\phi_\mu^2  &  \mp \phi_\mu^2   \\
\pm \phi_\mu^2           &  \mp \mu \pm 2\phi_\mu^2   
\end{array} \right) 
=: \mathcal{D} + \widetilde{\Mop}_\mu .
\end{equation}
Note that  $J^* = 2J^{-1}$ and so
\[ \text{Sp}(\Lop_\mu) = i\,\text{Sp}(\Mop_\mu), \quad \forall \mu \in [\mp \pi^2,+\infty). \]

\subsection{Basic spectral properties}

In this section, we recall basic spectral properties of the operators $\Lop_\mu$ and $\Mop_\mu$.
For this article to be self-contained, we propose proofs in Appendix B.

\begin{prop} \label{prop:basic_spect}
Let $\mu \in[\mp \pi^2,+\infty)$.
\begin{enumerate}

\item The spectrum of $\Mop_\mu$ and $\Mop_\mu^*$ is purely discrete and
the systems of eigenvectors and generalized eigenvectors for $\Mop_\mu$ and $\Mop_\mu^*$ 
(and hence for $\Lop_\mu$ and $\Lop_\mu^*$) form Schauder bases for $L^2((0,1),\C ^2)$.

\item All non-zero eigenvalues of $\Lop_\mu$ are purely imaginary:
$\text{Sp}(\Lop_\mu)=\{ \pm i \beta_{n,\mu} ; n \in \mathbb{N} \}$ where  
$(\beta_{n,\mu})_{n \in \mathbb{N}} \subset [0,+\infty)^\mathbb{N}$ is non decreasing
(here, multiple eigenvalues are repeated).

\item There exists $n_*=n_*(\mu) \in \mathbb{N}$ and $C=C(\mu)>0$ such that
\begin{equation} \label{Asymp_vap_t1}
|\beta_{n,\mu} - (n+n^*)^2\pi^2 | \le C  , \quad \forall n \in \mathbb{N}.
\end{equation}

\item The function $\mu \mapsto \beta_{n,\mu}$ is continuous for every $n \in \mathbb{N}$ and
\begin{equation} \label{vap_cv_-pi2}
\beta_{n,\mp \pi^2} = [(n+1)^2-1] \pi^2, \quad \forall n \in \mathbb{N}.
\end{equation}

\item The multiplicity of the eigenvalues of $\Lop_\mu$ is at most two. 
No non-zero eigenvalue possesses a generalized eigenvector. 

\item The vectors
\[ \Zn_0^+ = {0\choose \phi_\mu } \quad \text{ and } \quad  \Zn_0^- = {\partial _\mu \phi_\mu \choose 0} \] 
satisfy
\begin{equation} \label{jordan}
\Lop_\mu  \Zn_0^-=\Zn_0^+, \qquad \qquad  \Lop_\mu  \Zn _0^+ =0.
\end{equation}
Moreover $(\Zn_0^+,\Zn_0^-)$ is a basis of the generalized null space for $\Lop_\mu$. The vectors
\begin{equation} \label{def:Psi0+}
\Psi _0^- = {\phi_\mu \choose 0 }, \quad \text{ and } \quad  \Psi _0^+ = { 0 \choose \partial _\mu \phi_\mu } 
\end{equation}
satisfy
\begin{equation} \label{jordan*}
\Lop_\mu^*\Psi _0^+ = \Psi _0^-, \qquad \qquad \Lop_\mu^* \Psi _0^- =0.
\end{equation}
Moreover, $(\Psi_0^+,\Psi_0^-)$ is a basis of the generalized null space of $\Lop_{\mu}^*$.

\item Let $(\Zn_n^+)_{n \in \mathbb{N}^*}$ be normalized (see remark \ref{rk:norm} below)
eigenvectors of $\Lop_\mu$ associated to the eigenvalues $(+i\beta_{n,\mu})_{n \in \mathbb{N}^*}$
and $\Zn_n^-:=\overline{\Zn_n^+}$, then
$$\Lop_\mu \Zn _n^\pm = \pm i\beta _{n,\mu} \Zn_n^\pm  , \quad \forall n \in \mathbb{N}^*.$$
Let $(\Psi_n^+)_{n \in \mathbb{N}^*}$ be normalized eigenvectors of $\Lop_\mu^*$ associated to the eigenvalues $(-i\beta_{n,\mu})_{n \in \mathbb{N}^*}$
and $\Psi_n^-:=\overline{\Psi_n^+}$, then
$$\Lop_\mu^* \Psi _n^\pm = \mp i\beta _{n,\mu} \Psi_n^\pm  , \quad \forall n \in \mathbb{N}^*.$$
Moreover, if all non zero eigenvalue of $\Lop_\mu$ is simple then
\begin{equation} \label{Kronecker}
\langle \Zn_m^\sigma ,\Psi _n^\tau \rangle = 
\delta _{m,n}^{\sigma , \tau} := 
\left\{ \begin{array}{cl} 1, & m=n \, \, \textrm{and} \, \, \sigma = \tau, \\
 0, & \textrm{otherwise,} \end{array} \right.
\qquad  \forall m,n\in \mathbb{N}, \sigma,\tau \in \{ +,-\}   
\end{equation}
where the inner product is defined by (\ref{bracket}).

\item Let
$$ V_n^\mp := J \Phi_n^\pm,  \qquad \qquad  W_n^\mp:=J \Psi_n^\pm, \quad \forall n \in \mathbb{N}^*,$$
then,
$$\Mop_\mu   V_n^\pm = \pm \beta_{n,\mu} V_n^\pm , \qquad \qquad  \Mop_\mu^* W_n^\pm \, = \pm \beta_{n,\mu} W_n^\pm, \quad \forall n \in \mathbb{N}^*.$$ 
\end{enumerate}
\end{prop}

\begin{rk} \label{rk:exception}
When we use the vectors $\Zn_{n}^\pm$, $\Psi_{n}^\pm$, $V_{n}^\pm$, $W_{n}^\pm$
the symbols '$\pm$' and '$\mp$' do not refer to a distinction between the focusing and defocussing cases,
but to the sign of the associated eigenvalue.
\end{rk}

\begin{rk} \label{rk:norm}
Note that in the previous statement, the vectors $\Zn_n^\sigma$ and $\Psi _n^\sigma$ are defined up to 
a constant $c_n^\sigma \neq 0$, for every $n \geqslant 1$ and $\sigma \in \{\pm\}$:
$\langle c_m^\sigma \Zn_m^\sigma , \Psi _n^\tau / c_n^\tau \rangle   = \delta _{m,n}^{\sigma , \tau}$
for every sequence $(c_n^\pm)_{n \in \mathbb{N}^*} \subset \mathbb{R}^*$.
The 'normalization' we refer to in the statement \emph{(v)} will be chosen in Proposition \ref{Lkam}.
\end{rk}

\begin{rk}
We should have written  $\Zn_{n,\mu}^\pm$, $\Psi_{n,\mu}^\pm$, $V_{n,\mu}^\pm$, $W_{n,\mu}^\pm$ because these
vectors depend on $\mu$. We do not precise $\mu$ in subscript in order to simplify the notations.
\end{rk}

\subsection{Asymptotics of eigenvectors}

In the sequel, we use the $\mathcal{O} $-notation for uniform estimates:
\[ f_n(x) = g_n(x) +\mathcal{O} (n^{-\alpha }) \]
is to mean that there exists a constant $C$ and functions $R_n(x)$ such that
 \[ f_n(x) = g_n(x) + R_n(x) n^{-\alpha } \quad \textrm{and} \quad 
 |R_n(x)| \le C, \forall x \in (0,1), \forall n \in \mathbb{N}.\]

\begin{prop} \label{Lkam}
Let $\mu \in (\mp \pi^2,+\infty)$ and $n_*=n_*(\mu) \in \mathbb{N}$ be as in (\ref{Asymp_vap_t1}).
The normalization of $(\Zn_n^\pm,\Psi_n^\pm)$ may be chosen such that 
\begin{subequations}
\begin{eqnarray}
 V_n^\pm (x) &=& \st \sin [(n+n_*)\pi x] e^\pm  +\mathcal{O} \left( 1/n \right), \label{Wna} \\ 
 W_n^\pm (x) &=& \st \sin [(n+n_*)\pi x] e^\pm  +\mathcal{O} \left( 1/n \right), \label{Wnb}
\end{eqnarray} 
\end{subequations}
where $e^+ = \textstyle {1\choose 0}$ and $e^-= \textstyle {0\choose 1}$.
\end{prop}

\noindent \textbf{Proof of Proposition \ref{Lkam}:} 
In this proof, we omit the $\mu$ in subscripts, in order to simplify the notations,
and we deal with the focusing case (the defocussing case may be treated similarly).
First, we prove the estimate on 
\begin{equation} \label{def:uv}
V_n^+(x):={u_n(x) \choose v_n(x) }.
\end{equation}
The equation $\Mop_\mu V_n^+ = \beta_n V_n^+$ gives
\begin{subequations}
\begin{eqnarray}
 u_n'' +(\beta_n - \muz) u_n &=& -\phi_\mu^2 (2u_n+v_n), \qquad  u_n(0)=u_n(1)=0,  \\
 v_n'' -(\beta_n + \muz) v_n &=& -\phi_\mu^2 (u_n+2v_n), \qquad  v_n(0)=v_n(1) =0 \label{b}.
\end{eqnarray}   
\end{subequations}
For $n$ large enough, $(\beta_n+\mu)$ is positive (see (\ref{Asymp_vap_t1})), thus $\omega _n :=\sqrt{\beta _n+\muz }$ is well defined.  
From the relations
$$\left\lbrace \begin{array}{l}
u_n'' + [(n+n_*)\pi]^2 u_n = f_n(x) := \Big([(n+n_*)\pi]^2 - \beta_n + \mu - 2 \phi_\mu^2 \Big) u_n -  \phi_\mu v_n, \\
u_n(0)=u_n(1)=0
\end{array}\right.$$
we deduce that
\begin{equation} \label{u:MVC}
u_n(x) =  c \sin [(n+n_*)\pi x] +\frac{1}{(n+n_*)\pi } \int _0^x \sin [(n+n_*)\pi (x-\sigma )] f_n(\sigma ) d\sigma
\end{equation}
for some constant $c \in \mathbb{R}$ that may be taken equal to $2$ (see Remark \ref{rk:norm}).
We deduce from (\ref{Asymp_vap_t1}) that $u_n(x)=2\sin[(n+n_*)\pi x] + \mathcal{O}(1/n)$.
From (\ref{b}), we deduce that
\begin{equation} \label{v_Green}
v_n(x) = - \int _0^1 G_{\omega _n}  (x,\sigma) \phi_\mu(\sigma)^2 [u_n(\sigma)+2v_n(\sigma)] d\sigma,
\end{equation}
where
$$G_{\omega _n}  (x,\sigma)=-\frac{\sinh(\omega_n x) \sinh[\omega_n(1-\sigma)]}{\omega_n \sinh(\omega_n)} + \frac{\sinh[\omega_n(x-\sigma)]}{\omega_n} 1_{\sigma<x}.$$
The function $|G_{\omega _n } (x,\sigma)|$ assumes its maximum on $[0,1]^2$ at the point 
$(x,\sigma) = \left(\frac{1}{2} , \frac{1}{2} \right) $ and its  maximum value is given by 
\begin{equation} \label{omega}
 \left|G_{\omega _n } \left({\textstyle\frac{1}{2},\frac{1}{2} }\right)\right| =
\frac{\sinh ^2(\frac{\omega _n  }{2})}{\omega _n  \sinh(\omega _n  )} = \frac{\cosh(\omega _n  )-1}{2\omega _n  \sinh(\omega _n  )} = \mathcal{O}\left( 1/\omega_n \right).
\end{equation}
Thus, (\ref{v_Green}) and (\ref{Asymp_vap_t1}) justify that $v_n(x)=\mathcal{O}(1/n)$.
\\

The estimate on $V_n^-$ follows because
\begin{equation}
 V_n^- = {\overline{v_n} \choose \overline{u_n}}. 
\end{equation} 

Working similarly, we get the existence of a constant $C_n$ such that
$$ W_n^\pm (x) = \st C_n \sin[(n+n_*)\pi x] e^\pm +\mathcal{O} \left( 1/n \right).$$
Thus,
$$\begin{array}{ll}
\delta_{\sigma,\tau}^{n,m} 
& = \langle \Phi_m^\sigma , \Psi_n^\tau \rangle \\
& = \langle J^{-1} V_m^{\sigma'} , J^{-1} W_n^{\tau'} \rangle \\
& = \frac{1}{2} \langle  V_m^{\sigma'} , W_n^{\tau'} \rangle \\
& = \frac{1}{2} \langle 2 \sin[(m+n_*)\pi x] e^{\sigma'} + O\left(\frac{1}{m}\right) , 2 C_n \sin[(n+n_*)\pi x] e^{\tau'} + O \left( \frac{1}{n} \right) \rangle \\
& = \delta_{\sigma,\tau} 2 C_n \int_0^1 \sin[(m+n_*)\pi x] \sin[(n+n_*) \pi x] dx + O\left( \frac{1}{m}+\frac{1}{n} \right) \\
& = C_n \delta_{\sigma,\tau}^{n,m} + O\left( \frac{1}{m}+\frac{1}{n} \right).
\end{array}$$
Thus  $C_n=1+O(1/n)$ when $n \rightarrow + \infty$, which gives the conclusion.    \hfill $\Box $

\begin{prop} \label{Pr4} 
Let $\mu \in (\mp \pi^2,+\infty)$ and $n_*=n_*(\mu) \in \mathbb{N}$ be as in (\ref{Asymp_vap_t1}).
We denote $V_n^+(x) = {u_n(x)\choose v_n(x) }$ and $W_n^+(x) = {w_n(x)\choose z_n(x) }$.
There exist $\rho_n, \sigma_n, \widetilde{\rho}_n, \widetilde{\sigma_n} \in C^1([0,1],\mathbb{C})$, and $C>0$ such that  
\begin{equation} \label{bound:rho_r}
\|\rho_n\|_{C^1([0,1])}, \|\widetilde{\rho}_n\|_{C^1([0,1])}, 
\|\sigma_n\|_{C^1([0,1])}, \|\widetilde{\sigma}_n\|_{C^1([0,1])} \le C, \quad \forall n \in \mathbb{N}^*,
\end{equation}
\begin{equation} \label{rho_01}
\rho_n(0)=\rho_n(1)=\widetilde{\rho}_n(0)=\widetilde{\rho}_n(1)=0, \quad \forall n \in \mathbb{N}^*,
\end{equation}
\begin{equation} \label{u1} 
u_n(x) =  2\sin [(n+n_*)\pi x] +\frac{\sin [(n+n_*)\pi x]}{(n+n_*)\pi } \rho _n(x) - \frac{\cos [(n+n_*)\pi ]}{(n+n_*)\pi } \sigma _n(x) + \mathcal{O} \left( 1/n^2 \right), 
\end{equation}
\begin{equation}  \label{w1}
w_n(x) =  2\sin [(n+n_*)\pi x] +\frac{\sin [(n+n*)\pi x]}{(n+n_*)\pi } \widetilde{\rho}_n(x) - \frac{\cos [(n+n_*)\pi x]}{(n+n_*)\pi } \widetilde{\sigma}_n(x) 
+ \mathcal{O}\left( 1/n^2 \right),
\end{equation} 
\begin{equation} \label{v1}
v_n(x), z_n(x) = \mathcal{O}\left( 1/n^2 \right).
\end{equation}
\end{prop} 

\noindent \textbf{Proof of Proposition \ref{Pr4}:} In this proof, we omit $\mu$ in subscripts to simplify the notations
and we deal with the focusing case (the defocussing case may be treated similarly).
From (\ref{u:MVC}) and (\ref{Wna}) we get
$$\begin{array}{l}
 u_n(x) =  2 \sin[(n+n_*)n\pi x] + \mathcal{O}(1/n^2)
\\ 
+ \frac{1}{(n+n_*)\pi} \int_0^x \sin[(n+n_*)\pi(x-s)] 
\{ [(n+n_*)\pi]^2 - \beta_n + \mu - \phi_\mu(s)^2 \} 2 \sin[(n+n_*)\pi s] ds.
\end{array}$$
By developing $\sin[(n+n_*)\pi(x-s)]$, we get the conclusion with
$$\begin{array}{l}
\rho_n(x) :=  2 \int _0^x \cos [(n+n_*)\pi s]\sin[(n+n_*)\pi s]  \{ [(n+n_*)\pi]^2 - \beta_n + \mu - \phi_\mu(s)^2 \} ds,
\\ 
\sigma_n(x) :=  2 \int _0^x \sin^2[(n+n_*)\pi s] \{ [(n+n_*)\pi]^2 - \beta_n + \mu - \phi_\mu(s)^2 \} ds,
\end{array}$$
that satisfy (\ref{bound:rho_r}) (see (\ref{Asymp_vap_t1})). Note that $\rho_n(1)=0$ as the integral of an odd function. 
The decomposition of $w_n$ may be proved similarly. 
From (\ref{Wna}), (\ref{v_Green}) and (\ref{omega}), we get
$$v_n(x)=-2 \int_0^1 G_{\omega_n}(x,\sigma) \phi_\mu^2(\sigma) \sin[(n+n_*)\pi\sigma] d\sigma + \mathcal{O}\left( 1/n^2 \right).$$
Developing the hyperbolic sinuses, we get
\begin{equation} \label{Green_expanded}
\begin{array}{ll}
G_{\omega_n}(x,\sigma)=
&
\frac{1}{2\omega_n} \left[ 
\Big( -e^{\omega_n(x-\sigma)} + e^{\omega_n(x+\sigma-2)} + e^{-\omega_n(x+\sigma)}  \Big)\left(1+\mathcal{O}\left( 1/n \right) \right) 1_{\sigma>x}
\right. \\ 
& \left.
\Big( e^{\omega_n(x+\sigma-2)} + e^{-\omega_n(\sigma+x)} - e^{-\omega_n(x-\sigma)}  \Big)\left(1+\mathcal{O}\left( 1/n \right) \right) 1_{\sigma<x}
\right].
\end{array}
\end{equation}
In particular, $v_n(x)$ contains terms of the form
$$\begin{array}{ll}
& \frac{1}{2\omega_n} \int_x^1 e^{(-\omega_n\pm i (n+n_*) \pi)\sigma} \phi_\mu^2(\sigma) d\sigma e^{\omega_n x}
\\ = & 
-\frac{1}{2\omega_n} \left( 
\frac{e^{\pm i (n+n_*) \pi x}}{-\omega_n \pm i (n+n_*) \pi}  \phi_\mu^2(x) 
+ \int_x^1 \frac{e^{\omega_n(x-\sigma) \pm i (n+n_*) \pi \sigma}}{-\omega_n \pm i (n+n_*) \pi} (\phi_\mu^2)'(\sigma) d\sigma       \right)
\\  = &  \mathcal{O} \left( 1/n^2 \right).
\end{array}$$
Working similarly on the other terms of the right hand side of (\ref{Green_expanded}), we get $v_n(x)=\mathcal{O}(1/n^2)$. 
The estimates on $w_n$ and $z_n$ may be proved similarly. $\hfill \Box$

\subsection{Link with $H^3_{(0)}(0,1)$}

\begin{prop} \label{Prop:H3}
Let $\mu \in (\mp \pi^2,+\infty)$. There exists $C=C(\mu)>0$ such that, 
\begin{equation} \label{Ineq:H3}
\left( \sum\limits_{n=1}^\infty | n^3 \langle Z , \Psi_n^\pm \rangle |^2 \right)^{1/2}
\leqslant C \|Z\|_{H^3_{(0)}}, \quad \forall Z \in H^3_{(0)}((0,1),\mathbb{C}^2).
\end{equation}
\end{prop}

\noindent \textbf{Proof of Proposition \ref{Prop:H3}:} 
In this proof, we omit $\mu$ in subscript to simplify the notations
and we deal with the focusing case (the defocussing case may be treated similarly).
Let $\mu \in (- \pi^2,+\infty)$.
\\

\noindent\emph{First step: Existence of $C>0$ such that}
$$\left( \sum\limits_{n=1}^\infty | n \langle Z , \Psi_n^\pm \rangle |^2 \right)^{1/2}
\leqslant C \|Z\|_{H^1_0}, \quad \forall Z \in H^1_0((0,1),\mathbb{C}^2).$$
For $ Z \in H^1_0((0,1),\mathbb{C}^2)$, we have
$\langle Z , \Psi_n^- \rangle = \langle \widetilde{Z} , W_n^+ \rangle$
where $\widetilde{Z}:=JZ/2 \in H^1_0((0,1),\mathbb{C}^2)$
by Proposition \ref{prop:basic_spect} \emph{(viii)}.
Using (\ref{v1}), we see that it is sufficient to prove that
$$\left( \sum\limits_{n=1}^\infty \left| n \int_0^1 f(x) w_n(x) dx \right|^2 \right)^{1/2}
\leqslant C \|f\|_{H^1_0}, \quad \forall f \in H^1_0((0,1),\mathbb{C}).$$
Using integrations by part, (\ref{w1}) and (\ref{bound:rho_r}), we get
$$\begin{array}{lll}
\left( \sum\limits_{n=1}^\infty \left| n \int_0^1 f(x) w_n(x) dx \right|^2 \right)^{1/2}
& \leqslant &
\left( \sum\limits_{n=1}^\infty \left| n \int_0^1 f(x) \sin[(n+n_*) \pi x] dx \right|^2 \right)^{1/2}
\\ & &
+ \left( \sum\limits_{n=1}^\infty \left| \int_0^1 f(x) \sin[(n+n_*) \pi x] \rho_n(x) dx \right|^2 \right)^{1/2}
\\ & &
+ \left( \sum\limits_{n=1}^\infty \left| \int_0^1 f(x) \cos[(n+n_*) \pi x] \sigma_n(x) dx \right|^2 \right)^{1/2}
\\ & \leqslant &
C \left( \sum\limits_{n=1}^\infty \left| \int_0^1 f'(x) \cos[(n+n_*) \pi x] dx \right|^2 \right)^{1/2}
\\ & &
+ C \left( \sum\limits_{n=1}^\infty \left| \frac{1}{n\pi} \int_0^1 (f\rho_n)'(x)(x) \cos[(n+n_*) \pi x] dx \right|^2 \right)^{1/2}
\\ & &
+ C \left( \sum\limits_{n=1}^\infty \left| \frac{1}{n\pi}  \int_0^1 (f \sigma_n)'(x) \sin[(n+n_*) \pi x] dx \right|^2 \right)^{1/2}.
\end{array}$$
Bessel-Parseval inequality gives the conclusion. 
\\

\noindent\emph{Second step: Proof of (\ref{Ineq:H3}).} For $Z \in H^3_{(0)}((0,1),\mathbb{C}^2)$ we have
$\langle Z , \Psi_n^- \rangle = i \langle \Lop_\mu Z , \Psi_n^- \rangle/ \beta_n$,
which gives the conclusion thanks to (\ref{Asymp_vap_t1}) and the first step. \hfill $\Box$

\subsection{Asymptotic estimates}

\begin{prop} \label{GammaProp}
For $\mu \in (\mp \pi^2,+\infty)$ and $n \in \mathbb{N}^*$ we define
\begin{equation} \label{def:Gamman}
\Gamma_{n,\mu}^+ :=  
\left\langle  
(x\phi_\mu)'
 \left( \begin{array}{r}
 1\\ 0 
\end{array} \right) ,
\Psi_{n,\mu}^-(x) \right\rangle 
= \int_0^1 (x\phi_\mu )'(x)  \overline{\Psi_{n,\mu}^{(1)}(x)} dx.
\end{equation}
For every $\mu \in (\mp \pi^2,+\infty)$, there exists $C=C(\mu)>0$ such that  
\begin{equation} \label{est}
\left| \Gamma _{n,\mu}^+ - \frac{(-1)^{n+n_*+1}\phi_\mu '(1)}{\pi n}  \right| \leqslant \frac{C}{n^2}, \quad \forall n \in \mathbb{N}^*
\end{equation}
where $n_*=n_*(\mu)$ is as in (\ref{Asymp_vap_t1}).
\end{prop} 

\begin{rk}
Note that $\phi_\mu'(1) \neq 0$; otherwise, $\phi_\mu'(0)$ would vanish (symmetry of $\phi_\mu$) 
and $\phi_\mu$ would be identically zero, because of the uniqueness in Cauchy-Lipschitz theorem.
Thus Proposition \ref{GammaProp} gives the asymptotic behavior: $\Gamma_{n,\mu}^+ \sim (-1)^{n+n_*+1}\phi_\mu '(1)/(\pi n)$ when $n \rightarrow + \infty$.
\end{rk}

\noindent \textbf{Proof of Proposition \ref{GammaProp}:} We deduce from (\ref{w1}) that
$$\begin{array}{ll}
\int_0^1 (x\phi_\mu)'(x) w_n(x) dx = 
& \int_0^1 (x\phi_\mu)'(x) 2\sin [(n+n_*)\pi x]  dx
\\ & + \int_0^1 (x\phi_\mu)'(x) 
\Big( \frac{\sin [(n+n_*)\pi x]}{(n+n_*)\pi} \widetilde{\rho}_n(x) - \frac{\cos [(n+n_*)\pi x]}{(n+n_*)\pi} \widetilde{\sigma}_n(x)  \Big) dx 
+ \mathcal{O}\left( \frac{1}{n^2} \right).
\end{array}$$
Integrating by part each of the 3 terms in the right hand side and using (\ref{bound:rho_r}) we get
\begin{equation} \label{int_wn}
\int_0^1 (x\phi_\mu)'(x) w_n(x) dx = \frac{(-1)^{n+n_*+1} 2 \phi_\mu '(1)}{\pi n}  + \mathcal{O} \left( 1/n^2 \right).
\end{equation}
Using (\ref{v1}) and Proposition \ref{prop:basic_spect}\emph{(viii)} we get the conclusion. \hfill $\Box$

\section{Controllability of the linearized system}
\label{sec:Cont_Lin}

\begin{prop} \label{Prop:Cont_Lin}
Let $\mu \in (\mp \pi^2,+\infty)$ be such that
\begin{itemize}
\item \textbf{(A)} all non zero eigenvalues of $\Lop_\mu$ are simple,
\item \textbf{(B)} $\Gamma_{n,\mu}^+ \neq 0, \forall n \in \mathbb{N}^*$ (see (\ref{def:Gamman}) for the definition of $\Gamma_{n,\mu}^+$).
\end{itemize}
Then the map $d\Theta_{T,\mu}(0):\dot{H}^1_0((0,T),\mathbb{R}) \rightarrow H^3_{(0)}(0,1) \cap (\phi_\mu e^{\pm i\mu T})^\perp$
has a continuous right inverse.
\end{prop}

Here, we use the notation
$$(\phi_\mu e^{\pm i\mu T})^\perp:=\left\{ \Psi \in L^2(0,1) ; \Re \left( \int_0^1 \overline{\Psi(x)} \phi_\mu(x) dx e^{\pm i\mu T} \right)=0 \right\}.$$

\noindent \textbf{Proof of Proposition \ref{Prop:Cont_Lin}:} By Proposition \ref{Prop:C1}, we have
$$d\Theta_{T,\mu}(0).U= \widetilde{\Psi}(T) e^{\pm i \mu T},$$
where $\widetilde{\Psi}$ solves (\ref{syst_lin_static}). Identifying $H^3_{(0)}((0,1),\mathbb{C})$ with $H^3_{(0)}((0,1),\mathbb{R}^2)$
(by decomposition in real and imaginary parts), we get
$$d\Theta_{T,\mu}(0).U= Z(T) e^{\pm i\mu T},$$
where $Z=(\Re \widetilde{\Psi},\Im \widetilde{\Psi}) \in C^0([0,T],H^3_{(0)}((0,1),\mathbb{R}^2)) \cap C^1([0,T],H^1_0((0,1),\mathbb{R}^2))$ solves (\ref{Z}).
\\
 
By Proposition \ref{prop:basic_spect} \emph{(i)} and \emph{(vii)}, we have
$$ Z(t) = c_{0}^{+}(t)\Zn_0^++c_{0}^{-}(t)\Zn_0^- +\sum _{n\in \mathbb{N}^*} [c_n(t)\Zn_n^+ +\overline{c_n(t)}{\Zn}_n^-] 
\quad \text{ in } L^2((0,1),\mathbb{C}^2), \quad \forall t \in [0,T],$$
where $c_0^\pm(t):=\langle Z(t),\Psi_0^\pm \rangle$ and
$c_n(t):=\langle Z(t),\Psi_n^+ \rangle \in C^1([0,T],\mathbb{C})$ for every $n \in \mathbb{N}$. 
From the equation (\ref{Z}) we deduce that
$$\begin{array}{l}
\dot{c} _{0}^{-}(t)  = U(t) \Gamma_{0,\mu}^-, \\
\dot{c} _{0}^{+}(t)  = c_{0}^{-}(t) + U(t)\Gamma _{0,\mu}^+,\\
\dot{c} _n(t)        = i\beta_{n,\mu} c_n(t) + U(t)\Gamma_{n,\mu}^+, \quad \forall n \in \mathbb{N}^*,
\end{array}$$
where $\Gamma_{0,\mu}^\pm:=\int_0^1 (x\phi_\mu)'(x) (\Psi_0^\pm)^{(1)}(x) dx$.
Solving these ODEs and using the assumption $\int_0^T U =0$, we get
$$\begin{array}{ll}
c_0^-(T)=0,\\
c_0^+(T)=\Gamma_{0,\mu}^- \int_0^T (T-t) U(t) dt,\\
c_n(T)=e^{i \beta_{n,\mu}T} \Gamma_{n,\mu}^+ \int_0^T U(t) e^{-i \beta_{n,\mu}t} dt ,  \quad \forall n \in \mathbb{N}^*.
\end{array}$$
Integrating by parts and using $U(0)=U(T)=0$ we get
$$\begin{array}{ll}
c_0^-(T)=0,\\
c_0^+(T)=\Gamma_{0,\mu}^- \int_0^T \frac{(T-t)^2}{2} \dot{U}(t) dt,\\
c_n(T)=e^{i \beta_{n,\mu}T} \frac{\Gamma_{n,\mu}^+}{i \beta_{n,\mu}} \int_0^T \dot{U}(t) e^{-i \beta_{n,\mu}t} dt ,  \quad \forall n \in \mathbb{N}^*.
\end{array}$$

By Proposition \ref{prop:moment_pb} in Appendix D and (\ref{Asymp_vap_t1}), there exists a continuous map
$L_T: \mathbb{R} \times l^2(\mathbb{N}^*,\mathbb{C})  \rightarrow  L^2((0,T),\mathbb{R})$ such that, for every
$(d_0,(d_n)_{n \in \mathbb{N}^*}) \in  \mathbb{R} \times l^2(\mathbb{N}^*,\mathbb{C})$,
the function $\nu:=L_T(d_0,(d_n))$ solves
$$\left\lbrace \begin{array}{l}
\int_0^T \nu(t)dt=\int_0^T (T-t) \nu(t) dt=0,\\
\int_0^T \frac{(T-t)^2}{2} \nu(t)dt=d_0,\\
\int_0^T \nu(t) e^{-i\beta_{n,\mu}t} dt = d_n, \forall n \in \mathbb{N}^*.
\end{array}\right.$$

For $\Psi_f \in  H^3_{(0)}(0,1)$ such that 
\begin{equation} \label{hyp:Psif}
\Re \left( \int_0^1 \overline{\Psi_f(x)} \phi_\mu(x) e^{\pm i\mu T} dx \right)=0,
\end{equation}  
we define  $d(\Psi_f):=(d_n)_{n \in \mathbb{N}}$ by
$$d_0:= \frac{\langle Z_f , \Psi_0^+ \rangle}{\Gamma_{0,\mu}^-},
\qquad \qquad
d_n:=\frac{i \beta_{n,\mu} \langle Z_f,\Psi_n^+ \rangle e^{-i \beta_{n,\mu} T}}{\Gamma_{n,\mu}^+}, \quad \forall n \in \mathbb{N}^*$$
where $Z_f:=(\Re[\Psi_f e^{\mp i\mu T}],\Im [\Psi_f e^{\mp i\mu T}])$.

We remark that $\Gamma_{0,\mu}^- \neq 0$; indeed the relation (\ref{def:Psi0+}) and integrations by parts justify that
$$\Gamma_{0,\mu}^- = \int_0^1 (x\phi_\mu)'(x) \phi_\mu(x) dx = \frac{1}{2} \int_0^1 \phi_\mu(x)^2 dx >0.$$
By assumption \textbf{\emph{(B)}}, $d_n$ is well defined for every $n \in \mathbb{N}^*$.
Using (\ref{def:Psi0+}) and (\ref{hyp:Psif}), we see that $d_0 \in \mathbb{R}$.
By Proposition \ref{Prop:H3} and (\ref{est}), the map $\Psi_f \in H^3_{(0)} \cap (\phi_\mu e^{\pm i\mu T})^\perp \mapsto d(\Psi_f)$ takes values in $l^2(\mathbb{N},\mathbb{C})$.
We get the conclusion with $d\Theta_{T,\mu}(0)^{-1}.\Psi_f := L_T[ d(\Psi_f) ]$. $\hfill \Box$

\section{Genericity}
\label{sec:Generic}

In this section we verify that the assumptions \emph{\textbf{(A)}} and \emph{\textbf{(B)}}  in Proposition \ref{Prop:Cont_Lin}
hold generically with respect to the parameter $\muz $.

\begin{prop} \label{Prop:Generic}
There exists a countable set $J \subset (\mp \pi^2,+\infty)$ such that,
for every $\mu \in (\mp \pi^2,+\infty) \setminus J$,  all non zero eigenvalues of $\Lop_\mu$ are simple,
and $\Gamma_{n,\mu}^+ \neq 0, \forall n \in \mathbb{N}^*$.
\end{prop}

\subsection{Reformulation of the problem}

The purpose of the next two statements is to recast conditions  \emph{\textbf{(A)}} and \emph{\textbf{(B)}} such  that they become amenable to 
complex-variable methods. This is accomplished in (\ref{rk:Prop:SC}) below.

\begin{prop} \label{prop:Gamma_explicit}
Let $\mu \in (\mp \pi^2,+\infty)$.
We denote $\Psi_n^\pm={ f_{n,\mu}(x) \choose \mp i g_{n,\mu}(x) }$.
Then
$$\Gamma_{n,\mu}^+ = \frac{\phi_\mu'(1) g_{n,\mu}'(1)}{\beta_{n,\mu}}, \quad \forall n \in \mathbb{N}^*.$$
\end{prop}

\noindent \textbf{Proof of Proposition \ref{prop:Gamma_explicit}:} 
In this proof, we omit $\mu$ in subscript to simplify the notations, and we treat the focusing case (the defocussing one may be treated similarly).
From the relation $\Lop_\mu \Psi_n^\pm = \mp i \beta_{n} \Psi_n^\pm$, we get
\begin{equation} \label{syst_vep}
\begin{array}{l}
f_n''+ (\phi_\mu^2 -\mu) f_n = \beta_{n,\mu} g_n, \qquad f_n(0)=f_n(1)=0, \\
g_n''+ (3\phi_\mu^2-\mu) g_n = \beta_{n,\mu} f_n, \qquad g_n(0)=g_n(1)=0.
\end{array}   
\end{equation}
So, integration by parts gives
\begin{eqnarray*}
\Gamma_n^+   & = &\int_0^1 (x\phi_\mu)(x) f_n(x) dx  \\
             & = &\frac{1}{\beta_{n}} \int_0^1 \Big([\partial_x^2 + 3\phi_\mu^2-\mu ]g_n\Big)(x) (x\phi_\mu)'(x) dx \\
             & = &\frac{\phi_\mu'(1) g_n'(1)}{\beta_{n}} + \frac{1}{\beta_{n}} \int_0^1 \Big([\partial_x^2 + 3\phi_\mu^2-\mu ](x\phi_\mu)'\Big)(x) g_n (x) dx.
\end{eqnarray*}
Moreover, using (\ref{def:phimu}), we get
$$[\partial_x^2 + 3\phi_\mu^2-\mu ](x\phi_\mu)' 
= x (\phi_\mu''+\phi_\mu^3-\mu \phi_\mu)' + 3 (\phi_\mu''+\phi_\mu^3-\mu \phi_\mu) + 2 \mu \phi_\mu 
=  2 \mu \phi_\mu$$
and
\begin{eqnarray*}
\int_0^1 \phi_\mu(x) g_n(x) dx
& = &\frac{1}{\beta_{n}} \int_0^1 [\partial_x^2 + \phi_\mu^2-\mu ] f_n \phi_\mu dx=  \frac{1}{\beta_{n}} \int_0^1 [\partial_x^2 + \phi_\mu^2-\mu ] \phi_\mu f_n  dx =0,
\end{eqnarray*}
which gives the conclusion. $\hfill \Box$

\begin{prop} \label{Prop:SC}
Let $\mu \in (\mp \pi^2,+\infty)$, $n \in \mathbb{N}^*$, $(f_{n,\mu}^{[1]},g_{n,\mu}^{[1]})$, $(f_{n,\mu}^{[2]},g_{n,\mu}^{[2]})$ be the solutions of
\begin{equation} \label{syst:fg}
\left\lbrace \begin{array}{l}
f'' \pm (\phi_\mu^2 -\mu) f = \beta_{n,\mu} g,\\
g'' \pm (3\phi_\mu^2-\mu) g = \beta_{n,\mu} f, 
\end{array}\right.
\end{equation}
associated to the following initial conditions at $x=0$,
\begin{equation} \label{IC:fg}
\begin{array}{l}
f_{n,\mu}^{[1]}(0)=g_{n,\mu}^{[1]}(0)=(g_{n,\mu}^{[1]})'(0)=0, \qquad (f_{n,\mu}^{[1]})'(0)=1, \\
f_{n,\mu}^{[2]}(0)=g_{n,\mu}^{[2]}(0)=(f_{n,\mu}^{[2]})'(0)=0, \qquad (g_{n,\mu}^{[2]})'(0)=1,
\end{array}
\end{equation}
and
$$A_{n,\mu}:=\left(\begin{array}{cc}
f_{n,\mu}^{[1]}(1) & f_{n,\mu}^{[2]}(1) \\
g_{n,\mu}^{[1]}(1) & g_{n,\mu}^{[2]}(1)
\end{array}\right).$$
\begin{enumerate}
\item If $A_{n,\mu}$ is not the zero matrix, then $i \beta_{n,\mu}$ is a simple eigenvalue of $\Lop_\mu$.
\item If the first column of $A_{n,\mu}$ is not the zero vector ${0 \choose 0 }$ then $\Gamma_{n,\mu}^+ \neq 0$.
\end{enumerate}
In particular,
\begin{equation} \label{rk:Prop:SC}
f_{n,\mu}^{[1]}(1) \neq 0 \quad \Rightarrow \quad \beta_{n,\mu} \text{ is simple and } \Gamma_{n,\mu}^+ \neq 0.
\end{equation}
\end{prop}

\noindent \textbf{Proof of Proposition \ref{Prop:SC}:}
Let  $(f_{n,\mu}^{[3]},g_{n,\mu}^{[3]})$, $(f_{n,\mu}^{[4]},g_{n,\mu}^{[4]})$ be the solutions of (\ref{syst:fg}) such that
$$\begin{array}{l}
g_{n,\mu}^{[3]}(0)=(f_{n,\mu}^{[3]})'(0)=(g_{n,\mu}^{[3]})'(0)=0, \qquad f_{n,\mu}^{[3]}(0)=1, \\
f_{n,\mu}^{[4]}(0)=(f_{n,\mu}^{[4]})'(0)=(g_{n,\mu}^{[4]})'(0)=0, \qquad g_{n,\mu}^{[4]}(0)=1.
\end{array}$$
We assume that $i \beta_{n,\mu}$ is not a simple eigenvalue of $\Lop_\mu^*$.
Then (see Proposition \ref{prop:basic_spect} \emph{(v)})
there exists two linearly independent solutions 
$(f_{n,\mu}, g_{n,\mu})$ and $(\tilde{f}_{n,\mu},\tilde{g}_{n,\mu})$ of (\ref{syst_vep}).
They may be expanded with respect to the fundamental system
$${ f_{n,\mu} \choose g_{n,\mu}} = \sum_{k=1}^4 a_k { f_{n,\mu}^{[k]} \choose g_{n,\mu}^{[k]} }
\text{ and }
{ \tilde{f}_{n,\mu} \choose \tilde{g}_{n,\mu}} = \sum_{k=1}^4 \tilde{a}_k { f_{n,\mu}^{[k]} \choose g_{n,\mu}^{[k]} }
\text{ with } a_k, \tilde{a}_k \in \mathbb{C} \text{ for } k=1,...,4.$$
From the property $f_{n,\mu}(0)=g_{n,\mu}(0)=\tilde{f}_{n,\mu}(0)=\tilde{g}_{n,\mu}(0)=0$ we deduce that $a_3=a_4=\tilde{a}_3=\tilde{a}_4=0$.
From the property $f_{n,\mu}(1)=g_{n,\mu}(1)=\tilde{f}_{n,\mu}(1)=\tilde{g}_{n,\mu}(1)=0$, we deduce that the two linearly independent vectors
${a_3 \choose a_4}$ and ${\tilde{a}_3 \choose \tilde{a}_4}$ belong to the kernel of $A_{n,\mu}$. Thus $A_{n,\mu}=0$.
This proves the first statement.
\\

We assume that $A_{n,\mu} \neq 0$ and $\Gamma_{n,\mu}^+=0$.
By Proposition \ref{prop:Gamma_explicit}, we have $g_{n,\mu}'(1)=0$.
Since the eigenvalue $i \beta_{n,\mu}$ is simple,
$\Psi_n^{+}$ is either odd or even. In particular $g_{n,\mu}'(0)=\pm g_{n,\mu}'(1)=0$.
Thus, $\Psi_n$ is collinear to ${f_{n,\mu}^{[1]} \choose g_{n,\mu}^{[1]}}$.
We deduce from the relation $\Psi_n^+(1)=0$ that the first column of $A_{n,\mu}$ vanishes.
This proves the second statement. $\hfill \Box$

\subsection{Analyticity of the eigenvalues}

In this section we state the analytic dependence of
the eigenvalues of $\Lop_\mu$ with respect to $\mu$.
Because of the non--selfadjointness of $\Lop_\mu$ and $\Mop_\mu$, 
this property is not at all obvious.
In fact, there are simple examples of analytic families of $2\times 2$ matrices whose eigenvalues are not analytic functions of the 
parameter 
(see e.g. \cite{hryniv1999}).
We therefore provide a proof in Appendix C.

\begin{prop} \label{Prop:vap_analytic}
There exists continuous functions $F_n:[\mp \pi^2,\infty) \rightarrow \mathbb{R}^*_+$, for $n\in \mathbb{N}^*$,
that are analytic on $(\mp \pi^2,\infty)$ such that
\begin{itemize}
\item $\{F_n(\mu); n \in \mathbb{N}^*\}=\{ \beta_{n,\mu} ; n \in \mathbb{N}^* \}$, for every $\mu \in (\mp \pi^2,\infty)$,
\item $F_n(\mp \pi^2) = [(n+1)^2-1]\pi^2$.
\end{itemize}
\end{prop}

\subsection{Proof of Proposition \ref{Prop:Generic}}

The proof of Proposition \ref{Prop:Generic} follows from (\ref{rk:Prop:SC}) and the next result.

\begin{prop} \label{Prop:generic_interm}
There exists a countable set $J \subset (\mp \pi^2,+\infty)$ such that,
for every $\mu \in (\mp \pi^2,+\infty) \setminus J$ and $n \in \mathbb{N}^*$, the solution
$(f_{n,\mu}^{[1]},g_{n,\mu}^{[1]})$ of (\ref{syst:fg})-(\ref{IC:fg}) satisfies $f_{n,\mu}^{[1]}(1) \neq 0$.
\end{prop}

\noindent \textbf{Proof of Proposition \ref{Prop:generic_interm}:} 
We treat the focusing case (the defocussing one may be treated similarly).
Let $(F_n)_{n \in \mathbb{N}^*}$ be as in Proposition \ref{Prop:vap_analytic}.
We denote by $(k_{n,\mu},h_{n,\mu})$ the solution of
\begin{equation} \label{eq:kh}
\left\lbrace \begin{array}{l}
k_{n,\mu}''+ (\phi_\mu^2 -\mu) k_{n,\mu} = F_n(\mu) h_{n,\mu},\\
h_{n,\mu}''+ (3\phi_\mu^2-\mu) h_{n,\mu} = F_n(\mu) k_{n,\mu},\\
k_{n,\mu}(0)=h_{n,\mu}(0)=h_{n,\mu}'(0)=0, \quad k_{n,\mu}'(0)=1,
\end{array}\right.
\end{equation}
and we introduce the map
$$\begin{array}{|cccc}
G_n: & [-\pi^2,+\infty) & \rightarrow & \mathbb{R} \\
     &     \mu          & \mapsto     & k_{n,\mu}(1).
\end{array}$$

\noindent \emph{First step: $G_n$ is analytic for every $n \in \mathbb{N}^*$.}
Let $n \in \mathbb{N}^*$. Since $n$ is fixed in all this step, we will write 
$k_\mu$, $h_\mu$, $F$ instead of $k_{n,\mu}$, $h_{n,\mu}$, $F_n$.
Let $\mu_0 \in (-\pi^2,+\infty)$. The functions $\mu  \mapsto \phi_\mu$ and $\mu \mapsto F(\mu)$
may be extended as holomorphic functions of $\mu \in \Omega$ where
$\Omega:=\{ \mu \in \mathbb{C} ; \mu_0-\epsilon < \Re(\mu) < \mu_0+\epsilon, -\epsilon < \Im(\mu) < \epsilon\}$ for some $\epsilon>0$,
by the sum of the converging Taylor series at $\mu_0$.
For $\mu \in \Omega$, we introduce the notations
$$\begin{array}{c}
\mu_1:=\Re(\mu), \quad \mu_2:=\Im(\mu), \quad F^{(1)}(\mu):=\Re[F(\mu)], \quad F^{(2)}(\mu):=\Im[F(\mu)],\\ 
a_1(x):=\Re[\phi_\mu(x)^2-\mu], \quad a_2(x):=\Im[\phi_\mu(x)^2-\mu], \\
b_1(x):=\Re[3\phi_\mu(x)^2-\mu], \quad b_2(x):=\Im[3\phi_\mu(x)^2-\mu],\\
k_{\mu}^{(1)}(x):=\Re[k_{\mu}(x)], \quad k_{\mu}^{(2)}(x):=\Im[k_{\mu}(x)], \quad h_{\mu}^{(1)}(x):=\Re[h_{\mu}(x)], \quad h_{\mu}^{(2)}(x):=\Im[h_{\mu}(x)].
\end{array}$$
We deduce from (\ref{eq:kh}) that, for every $\mu \in \Omega$,
\begin{subequations} 
\begin{eqnarray}
(k_\mu^{(1)})''+ a_1 k_\mu^{(1)} - a_2 k_\mu^{(2)} & = & F^{(1)} h_\mu^{(1)} - F^{(2)} h_\mu^{(2)}, \label{eq:k1} \\
(k_\mu^{(2)})''+ a_1 k_\mu^{(2)} + a_2 k_\mu^{(1)} & = & F^{(1)} h_\mu^{(2)} + F^{(2)} h_\mu^{(1)}, \label{eq:k2} \\
(h_\mu^{(1)})''+ b_1 h_\mu^{(1)} - b_2 h_\mu^{(2)} & = & F^{(1)} k_\mu^{(1)} - F^{(2)} k_\mu^{(2)}, \label{eq:h1} \\
(h_\mu^{(2)})''+ b_1 h_\mu^{(2)} + b_2 h_\mu^{(1)} & = & F^{(1)} k_\mu^{(2)} + F^{(2)} k_\mu^{(1)}. \label{eq:h2} \\
\end{eqnarray}
\end{subequations}
In particular, for every  $(\mu_1,\mu_2) \in \widetilde{\Omega}:=(\mu_0-\epsilon,\mu_0+\epsilon) \times (-\epsilon,\epsilon)$, the function 
$Y_\mu:=(k_{\mu}^{(1)}, k_{\mu}^{(2)}, h_{\mu}^{(1)}, h_{\mu}^{(2)})$ solves an equation of the form
\begin{equation} \label{eq:Ymu}
\left\lbrace \begin{array}{l}
\frac{d^2 Y_\mu}{dx^2}=\mathcal{F}(x,Y_\mu,\mu_1,\mu_2) \\
Y_\mu(0)=(0,0,0,0),\\
Y_\mu'(0)=(1,0,0,0).
\end{array}\right.
\end{equation}
where the function $\mathcal{F}$ is of class $C^1$ with respect to $(x,Y,\mu_1,\mu_2) \in (0,1) \times \mathbb{R}^4 \times \widetilde{\Omega}$,
thus $Y_\mu$ has partial derivatives with respect to $\mu_1$ and $\mu_2$.
Now, we prove that they satisfy the Cauchy-Riemann relations, in order to get the holomorphy of $\mu \in \Omega \mapsto Y_\mu(1)$. 
We introduce the functions
$$K_{i,j}:=\frac{\partial k_\mu^{(i)}}{\partial \mu_j}, \quad H_{i,j}:=\frac{\partial h_\mu^{(i)}}{\partial \mu_j}, \quad \forall i,j \in \{1,2\}.$$
Computing 
$\partial_{\mu_1}(\ref{eq:k1})-\partial_{\mu_2}(\ref{eq:k2})$,
$\partial_{\mu_2}(\ref{eq:k1})+\partial_{\mu_1}(\ref{eq:k2})$,
$\partial_{\mu_1}(\ref{eq:h1})-\partial_{\mu_2}(\ref{eq:h2})$,
$\partial_{\mu_2}(\ref{eq:h1})+\partial_{\mu_1}(\ref{eq:h2})$,
and using the Cauchy-Riemann relations on $a_1$, $a_2$, $b_1$, $b_2$, $F^{(1)}$, $F^{(2)}$, we get
$$\begin{array}{l}
(K_{1,1}-K_{2,2})''+a_1(K_{1,1}-K_{2,2})-a_2(K_{2,1}+K_{1,2})=F^{(1)}(H_{1,1}-H_{2,2})-F^{(2)}(H_{1,2}+H_{2,1}),\\
(K_{1,2}+K_{2,1})''+a_1(K_{1,2}+K_{2,1})-a_2(K_{2,2}-K_{1,1})=F^{(1)}(H_{1,2}+H_{2,1})-F^{(2)}(H_{2,2}-H_{1,1}),\\
(H_{1,1}-H_{2,2})''+b_1(H_{1,1}-H_{2,2})-b_2(H_{2,1}+H_{1,2})=F^{(1)}(K_{1,1}-K_{2,2})-F^{(2)}(K_{1,2}+K_{2,1}),\\
(H_{1,2}+K_{2,1})''+b_1(H_{1,2}+H_{2,1})-b_2(H_{2,2}-H_{1,1})=F^{(1)}(K_{1,2}+K_{2,1})-F^{(2)}(K_{2,2}-K_{1,1}),\\
(K_{1,1}-K_{2,2},K_{1,2}+K_{2,1},H_{1,1}-H_{2,2},H_{1,2}+H_{2,1})(0)=(0,0,0,0),\\
(K_{1,1}-K_{2,2},K_{1,2}+K_{2,1},H_{1,1}-H_{2,2},H_{1,2}+H_{2,1})'(0)=(0,0,0,0).
\end{array}$$
The uniqueness of the solution of this linear system ensures that
$K_{1,1}-K_{2,2}=K_{1,2}+K_{2,1}=H_{1,1}-H_{2,2}=H_{1,2}+H_{2,1}=0$.
In particular $(K_{1,1}-K_{2,2})(1)=(K_{1,2}+K_{2,1})(1)$, which proves the holomorphy of $G_n$ on $\Omega$.
\\

\noindent \emph{Second step: $G_n(-\pi^2) \neq 0$, $\forall n \in \mathbb{N}^*$.}
Let $n \in \mathbb{N}^*$. Thanks to (\ref{limit:phi_mu}) and the second conclusion of Proposition \ref{Prop:vap_analytic}, 
we have $G_{n+1}(-\pi^2)=f(1)$ where $(f,g)$ is the solution of the Cauchy problem
\begin{equation} 
\left\lbrace \begin{array}{l}
f''+ \pi^2 f = (n^2-1)\pi^2 g,\\
g''+ \pi^2 g = (n^2-1)\pi^2 f,\\
f(0)=g(0)=g'(0)=0, \quad f'(0)=1. 
\end{array}\right.
\end{equation}
This system may be written
$$\left\lbrace \begin{array}{l}
\left( \frac{d^2}{dx^2} + \pi^2 \right)^2 f = (n^2-1)^2 \pi^4 f, \\
f(0)=f''(0)=0, \quad f'(0)=1, \quad f^{(3)}(0)=-\pi^2.
\end{array}\right.$$
Thus, $f$ may be computed explicitly. In particular,
$$f(1)=\left\lbrace \begin{array}{l}
-(1 + 2/(\pi-3)) \quad \text{ if } n=1,\\
3/4  \quad \text{ if } n=2,\\
\sinh(\sqrt{n^2-2} \pi)/(2\pi\sqrt{n^2-2})  \quad \text{ if } n \geqslant 3,
\end{array}\right.$$
which gives the conclusion.
\\

\noindent \emph{Third step: Conclusion.} 
By the isolated zero principle, for every $n \in \mathbb{N}^*$, there exists a countable set
$J_n \subset (-\pi^2,+\infty)$ such that,
$G_n(\mu) \neq 0$  for every $\mu \in (-\pi^2,+\infty) \setminus J_n$. 
Then, $J:=\cup_{n \in \mathbb{N}^*} J_n$ gives the conclusion. $\hfill \Box$

\section{Proof of the main result}
\label{sec:proof}

Let $J$ be as in Proposition \ref{Prop:Generic},  $\mu \in (\mp \pi^2,+\infty)\setminus J$ and $T>0$.
By Proposition \ref{Prop:C1}, the map 
$$\Theta_{T,\mu}:  \dot{H}^1_0((0,T),\mathbb{R})  \rightarrow  H^3_{(0)}(0,1) \cap \mathcal{S}_{\|\phi_\mu\|_{L^2}}$$
is $C^1$. By Proposition \ref{Prop:Cont_Lin} and (\ref{Prop:Generic}), 
$$d\Theta_{T,\mu}(0): \dot{H}^1_0((0,T),\mathbb{R})  \rightarrow  H^3_{(0)}(0,1) \cap (\phi_\mu e^{\pm i\mu T})^\perp$$
has a continuous right inverse. The inverse mapping theorem gives the conclusion.

\section{Conclusion and perspectives}
\label{sec:ccl}

Motivated by the control of Bose--Einstein condensates, we have studied the controllability  of the 
nonlinear Schr\"odinger equation (focusing and defocusing) 
with a bilinear control term arising from manipulating the size of a ``hard-wall" (box) trap. 
We showed that
local exact controllability around the ground state holds generically with respect to the parameter $\muz $. 
Since $\muz $ is a parameter associated with the 
\textit{transformed} problem (\ref{syst}), this leaves the question of whether genericity also holds with respect to the 
system parameter $\scat $ of the \textit{original} problem (\ref{GPE}). This is indeed so, as  is readily seen from the identity 
\begin{equation*}
  \| \phi _\muz \| _{L^2(0,1)}^2  
  =  \frac{2\scat m}{\hbar ^2} \label{bound4} 
   \end{equation*}  
and the convexity condition (\ref{convexite})\footnote{Another possible parameter is the initial and final size 
of the trap, which in (\ref{BC:L}) was set to one for convenience.}. 
While the genericity property implies that 
local controllability holds with ``probability one w.r.t. \textit{random} choices" of $\muz $ (or $\scat $),
 for any \textit{particular} value of  $\muz $  (resp. $\scat$) Theorem \ref{thm:Main} cannot be applied directly. 
 It will be shown elsewhere \cite{CT13} how rigorous numerical 
  computation can be utilized in these cases.

Of the numerous possible generalizations of the control problem considered in the present paper we briefly mention three:
\begin{enumerate}
\item more \textit{general nonlinearities};
\item controllability around \textit{excited} states\footnote{These are (real-valued) solutions of (\ref{def:phimu}), 
with a positive number of zeros (``nodes'') within the interval
$(0,1)$ \cite{CCR,CCR00}. (The node-less solution is the ground state.)};
 \item \textit{global} exact controllability. 
 \end{enumerate}
In (i) and (ii) several steps of our 
approach will need to be adapted, such as the study of the spectrum of the operator $\Lop_\mu$, 
which may no longer be purely imaginary, or the proof of the genericity result in Section \ref{sec:Generic}, which  
uses the convexity inequality (\ref{convexite}). We conjecture that (i) can be handled 
for ``benign'' cases such as certain power nonlinearities and that (ii) holds at least in the defocusing case. 
To prove (iii) one may try to adapt the techniques of \cite{Nersesyan2}, although, due to 
the nonlinearity of the equation, significant new ideas will be required.

\appendix

\section{Ground states: proof}
\label{appendix:GS}

In this section we prove Proposition \ref{Prop:GS}.

First, we treat the focusing case.
Let $\mu \in (-\pi^2,+\infty)$. There exists a unique solution $w_\mu \in M$ of the minimization problem
\begin{equation} \label{min_pb_focus}
\begin{array}{c}
J_\mu(w_\mu) = \inf\left\{ J_\mu(\varphi)  ; \varphi \in M \right\},\\
J_\mu(\varphi):=\int_0^1 [\varphi'(x)^2 + \mu \varphi(x)^2] dx, \\
M:=\left\{  \varphi \in H^1_0((0,1),\mathbb{R}); \int_0^1 \varphi(x)^4 dx =1 \right\},
\end{array}
\end{equation}
and a Lagrange multiplier $\alpha_\mu \in \mathbb{R}$ such that
\begin{equation} \label{EL_focus}
\left\lbrace \begin{array}{l}
-w_\mu'' + \mu w_\mu=\alpha_\mu w_\mu^3, \\
w_\mu(0)=w_\mu(1)=0.
\end{array}\right.
\end{equation}
Then
$$\alpha_\mu \int_0^1 w_\mu(x)^4 dx = \int_0^1 \Big( w_\mu'(x)^2 + \mu w_\mu(x)^2 \Big) dx >0;$$
thus $\alpha_\mu >0$ and $\phi_\mu := \sqrt{\alpha_\mu} \varphi_\mu$ gives the solution. 
An explicit formula of $\phi_\mu$ is available in terms of Jacobian elliptic functions.
For $\mu \in (-\pi^2,+\infty)$, we first find the solution $k=k(\mu)$ of the equation
$$\mu  = 4 (2k^2 - 1) K(k)^2$$
where $K$ denotes the complete elliptic integral of the first kind (see e.g. \cite{AS65}).
Note that the function $K:[0,1),\rightarrow [\pi/2,+\infty)$ is continuous, analytic on $(0,1)$, bijective and $K'>0$ on $(0,1)$.
Thus, the reciprocal $k=k(\mu)$ defines a function
$k:[-\pi^2,+\infty) \rightarrow [0,1)$ continuous, analytic on $(-\pi^2,+\infty)$, bijective with $k'>0$ on $(0,+\infty)$.
Then, the function $\phi_\mu$ is given by the formula \cite{CCR00}
$$\phi_\mu(x)=2 \sqrt{2} k K(k) \, \text{cn} \left( 2 K(k) \Big( x-\frac{1}{2} \Big) , k \right),$$
where $\text{cn}$ is the elliptic cosine function.
This proves the analyticity of the map $\mu \in (-\pi^2,+\infty) \mapsto \phi_\mu \in L^2(0,1)$ and the relation
$$ \int_0^1 \phi_\mu(x)^2 dx = 8 k^2 K(k) \int_{0}^{K(k)} \text{cn}(y)^2 dy = 8 K(k) F(k)$$
where $F(k):=E(k)-(1-k^2)K(k)$ and $E$ is the complete elliptic integral of the second kind.
The function $F$ is positive and satisfies $F'(k)=k K(k)$ on $(0,1)$, thus
$$2 \langle \partial_\mu \phi_\mu , \phi_\mu  \rangle
= 8 k'(\mu) \Big( K'[k(\mu)] F[k(\mu)] + K[k(\mu)] F'[k(\mu)] \Big)>0, \quad \forall \mu \in (-\pi^2,+\infty).$$
Note that, when $\mu$ tends to $-\pi^2$, then
$k(\mu) \rightarrow 0$, and $K[k(\mu)]$ is bounded, which proves (\ref{limit:phi_mu}). 
\\

In the defocussing case we will not need the variational description of the ground state, so we omit this point. 
Again, an explicit formula of $\phi_\mu$ is available in terms of Jacobian elliptic functions.
For $\mu \in (\pi^2,+\infty)$, we first find the solution $k=k(\mu)$ of the equation
$$\mu  = 4 (k^2 + 1) K(k)^2.$$
This defines a function $k:[\pi^2,+\infty) \rightarrow [0,1)$ continuous, analytic on $(\pi^2,+\infty)$, such that $k'>0$ on $(\pi^2,+\infty)$.
Then \cite{CCR} 
$$\phi_\mu(x):=2 \sqrt{2} k K(k)\, \text{sn} \Big( 2 K(k) x , k \Big),$$
where $\text{sn}$ is the Jacobian elliptic sine function.
This proves the analyticity of $\mu \in (\pi^2,+\infty) \mapsto \phi_\mu \in L^2(0,1)$ and the relation
$$\int_0^1 \phi_\mu(x)^2 dx = 4 k^2 K(k) \int_0^{K(k)} \text{sn}(y,k)^2 dy = 8 K(k) G(k)$$
where $G(k):=E(k)-K(k)$ is positive and satisfies $G'(k)=k E(k)/(1-k^2)$ for every $k \in (0,1)$.
The proof may be ended as above. $\hfill \Box$

\section{Basic spectral properties: proof}

In this appendix we provide the proof of Proposition \ref{prop:basic_spect}.
Our proof is similar to the one for the whole space case, which has been studied extensively; 
its elements are taken from \cite{K66}, \cite[Appendix B]{LT_MMAS} and adapted from \cite{RSS05}.

\subsection{Preliminaries}

\begin{prop} \label{Prop:techn}
In both focusing and defocussing cases, we have
\begin{equation} \label{Ker_L}
\text{Ker}(L_\mu^-)=\mathbb{C} \phi_\mu, \qquad
\text{Ker}(L_\mu^+)=\{0\}, \quad \forall \mu \in (\mp \pi^2,\infty).
\end{equation}
In the focusing case,
\begin{equation} \label{Ker_L+}
L_\mu^+ \text{ has only one negative eigenvalue, } \forall \mu \in (-\pi^2,\infty).  
\end{equation} 
In the defocussing case, 
\begin{equation} \label{Ker_L+_defocusing}
 L_\mu^+ > 0, \quad \forall \mu \in (\pi^2,+\infty).
\end{equation}
\end{prop}

\noindent \textbf{Proof of Proposition \ref{Prop:techn}:} 

\noindent \emph{First step: Proof of $\text{Ker}(L_\mu^-)=\mathbb{C} \phi_\mu$.} 
We recall that
$\text{Ker}[L_\mu^-]:=\{ w \in H^2 \cap H^1_0(0,1) ; (\partial_x^2 \pm \phi_\mu^2 \mp \mu)w \equiv 0 \}$.
The linear map
$$\begin{array}{|ccc}
\text{Ker}[L_\mu^-] & \rightarrow & \mathbb{C} \\
 w                  &  \mapsto    & w'(0)
\end{array}$$
is injective, thanks to the uniqueness in Cauchy-Lipschitz theorem.
Thus $\text{dim}[\text{Ker}(L_\mu^-)] \leqslant 1$. Clearly, $L_\mu^- \phi_\mu =0$, which gives the conclusion.
\\

\noindent \emph{Second step: Proof of $\text{Ker}(L_\mu^+)=\{0\}$ and (\ref{Ker_L+}) in the focusing case.}
This will be achieved thanks to Step 2.1, Step 2.2 and Step 2.3 below.
\\

\noindent \emph{Step 2.1: $L_\mu^+$ has at least one negative eigenvalue.}
This follows from 
$$\langle L_\mu^+\phi _\muz , \phi_\muz \rangle = - 2 \|\phi_\mu\|_{L^4}^4 <0$$
and the minimax principle.
\\

\noindent \emph{Step 2.2: $\langle L_\mu^+ \eta , \eta \rangle \geqslant 0, \forall \eta \perp \phi_\mu^3$.}
We use the characterization of $\phi_\mu$ by the minimization problem (\ref{min_pb_focus}).
Let $\eta \in H^2 \cap H^1_0((0,1),\mathbb{R})$ be such that $\eta \perp \phi_\mu^3$ in $L^2(0,1)$.
Let $z \mapsto w(.,z) \in  H^2 \cap H^1_0((0,1),\mathbb{R})$ be a smooth curve such that 
$w(.,0)=w_\mu$, 
$\dot{w}:=\partial_z [ w(.,z)]_{z=0}=\eta$,
$\| w(.,z)\|_{L^4(0,1)} \equiv 1$.
Since $w_\mu$ solves the minimization (\ref{min_pb_focus}), the function $z \mapsto J[w(.,z)]$ has its minimum at $z=0$; thus
$$0 = \frac{d}{dz} \Big[ J_\mu[w(.,z)] \Big]_{z=0} = \int_0^1 \Big( w_x \dot{w}_x + \mu w \dot{w} \Big) dx,$$
\begin{equation} \label{J_ineq}
0 \leqslant \frac{d^2}{dz^2} \Big[ J_\mu[w(.,z)] \Big]_{z=0} =  \int_0^1 \Big( \dot{w}_x^2 + w_x \Ddot{w}_x + \mu  \dot{w}^2 + \mu w \Ddot{w} \Big) dx.
\end{equation}
Moreover, $w(.,z) \in M$ for every $z$, thus
$$0 = \int_0^1 w^3 \dot{w} dx = \int_0^1  \Big( 3 w^2 \dot{w}^2 + w^2 \Ddot{w} \Big) dx.$$
The Euler-Lagrange equation (\ref{EL_focus}) and the previous relation give, at $z=0$,
$$\int_0^1  (-w_{xx}+ \mu w) \Ddot{w} dx = \alpha_\mu \int_0^1 w^3 \Ddot{w} dx = -3 \alpha_\mu \int_0^1 w^2 \dot{w}^2= - 3 \int_0^1 \phi_\mu^2 \eta^2.$$
Incorporating this relation in (\ref{J_ineq}) gives
\begin{eqnarray*}
0 & \leqslant & \int_0^1 \Big( (-\dot{w}_{xx}+\mu\dot{w})\dot{w} + (-w_{xx}+\mu w)\Ddot{w} \Big) dx \\
  & =  &\int_0^1 \Big( - \eta_{xx} + \mu \eta - 3 \phi_\mu^2 \eta \Big) \eta  dx\, = \,  \langle L_\mu^+ \eta , \eta \rangle.
\end{eqnarray*}

\noindent \emph{Step 2.3: The second eigenvalue $\lambda_2$ of $L_\mu^+$ is $>0$.}
Thanks to Step 2.2, we know that the second eigenvalue of $L_\mu^+$ is $\geqslant 0$.
Let us assume that $\lambda_2=0$. Let $v \in D(L_\mu^+)$ be such that $L_\mu^+ v=0$.
By symmetry of $\phi_\mu$ with respect to $x=1/2$, one may assume that $v$ is odd or even (with respect to $x=1/2$).
If $v$ is odd, then $v(1/2)=0$ and $v=c\phi_\mu$ for some $c \in \mathbb{R}$, thanks to ODE solutions uniqueness.
But this is impossible because $\phi_\mu'$ does not vanish at $x=0$ and $x=1$.
Thus $v$ is even. The function $v$ has one zero in $(0,1)$ (second eigenvalue of a Sturm-Liouville operator).
By symmetry, this must occur at $x=1/2$, which is impossible as saw before. 
\\

\noindent \emph{Third step: Proof of (\ref{Ker_L+_defocusing}) in the defocussing case.}
We prove by contradiction that the smallest eigenvalue is positive. To this end, let $E $
be the smallest eigenvalue, $u \in H^2 \cap H^1_0((0,1),\mathbb{R})\setminus \{0\}$ a corresponding eigenfunction, and assume $E\le 0$. Then $u$ may be assumed to be positive on $(0,1)$ because it is the ground state of $L_\mu^+$.
Thus $\langle u, \phi _\muz \rangle >0$ and so 
$$0\ge E \langle u, \phi _\muz \rangle = \langle L_\mu^+ u , \phi_\mu \rangle = \langle u , L_\mu^+ \phi_\mu \rangle = 2 \int_0^1 u(x) \phi_\mu(x)^3 dx >0,$$
which is impossible. Therefore $L_\mu^+>0$. $\hfill \Box$

\subsection{Statements \emph{(i)} and \emph{(iii)}}

The operator $\mathcal{D} $ defined by (\ref{DM})
is self-adjoint, with compact resolvent and simple eigenvalues with an infinite asymptotic gap:
$$\mathcal{D}  {\sin (n\pi x) \choose 0} = (n\pi)^2 {\sin (n\pi x) \choose 0}, \qquad
\mathcal{D}  {0 \choose \sin (n\pi x) } = -(n\pi)^2 {0 \choose \sin (n\pi x) }, \quad \forall n \in \mathbb{N}^*.$$
The operator $\widetilde{\mathcal{M}}_\mu$ is bounded on $L^2((0,1),\mathbb{C}^2)$.
By applying \cite[Chapter V, paragraph 3, Theorem 4.15.a on Page 293]{K66}),
we get the first statement of Proposition \ref{prop:basic_spect}
and the third one, assuming that the second one holds (which will be proved independently below).

\subsection{Statement \emph{(ii)}}

This proof follows the one of \cite[Lemma 12.11]{RSS05}, in the case of NLS on the whole line.
In order to simplify the notations, we do not write $\mu$ in subscript.
Let us consider the operator
\begin{equation} \label{Lop2}
\Lop^2=-\left( \begin{array}{cc}
\mathcal{T}^* & 0 \\
0             & \mathcal{T}
\end{array}\right)
\quad \text{  where  } \quad
D(\mathcal{T}):=H^4_{(0)}(0,1), \quad \quad \mathcal{T}:=L^+ L^-.
\end{equation}

\noindent \emph{First step: $\text{Sp}(\mathcal{T}) \subset \mathbb{R}$.}
Let $E \in \mathbb{C}-\{0\}$ be an eigenvalue of $\mathcal{T}$ and $\psi$ be an associated eigenvector: 
\begin{equation} \label{Tpsi=Epsi}
\mathcal{T}\psi=E\psi.
\end{equation}
Then, $\psi=\psi_1+c\phi_\mu$, where $\psi_1 \perp \phi_\mu$ and $\psi_1 \neq 0$ (because $L^-\phi_\mu=0$).
Thus, (\ref{Tpsi=Epsi}) gives
\begin{equation} \label{Tpsi=Epsi_bis}
\Big[(L^-)^{1/2} L^+ (L^-)^{1/2}\Big] \Big( (L^-)^{1/2} \psi_1 \Big) = E \Big( (L^-)^{1/2} \psi_1 \Big).
\end{equation}
Moreover $(L^-)^{1/2} \psi_1  \neq 0$ thanks to (\ref{Ker_L}).
Thus $E$ is an eigenvalue of the symmetric operator $(L^-)^{1/2} L^+ (L^-)^{1/2}$, so $E \in \mathbb{R}$.
\\

\noindent \emph{Second step: (\ref{Ker_L+}) implies $\text{Sp}(\mathcal{T}) \subset \mathbb{R}^+$ in the focusing case.}
The map
$$g(E):=\langle (L^+-E)^{-1} \phi_\mu , \phi_\mu \rangle$$
is well defined for $E \in (-E^*,0]$, where $-E^*$ is the negative eigenvalue of $L^+$.
Moreover, we have
$$g'(E)= \|(L^+-E)^{-1} \phi_\mu\|^2 \geqslant 0, \quad \forall E \in (-E^*,0),$$
$$g(0)=-\langle \partial_\mu \phi_\mu , \phi_\mu \rangle <0,$$
thanks to (\ref{convexite}), thus 
\begin{equation} \label{absurd}
g(E)<0, \quad \forall E \in (-E^*,0).
\end{equation}

We prove by contradiction that the eigenvalues of $\mathcal{T}$ are $\geqslant 0$.
We assume that $\mathcal{T}$ has a negative eigenvalue $E<0$.
From (\ref{Tpsi=Epsi_bis}), we deduce that $(L^-)^{1/2} L^+ (L^-)^{1/2}$ has a negative eigenvalue in $\text{Ker}[L^-]^{\perp}$;
there exists $\zeta \in \text{Ker}[L^-]^{\perp}$ such that
$$\langle (L^-)^{1/2} L^+ (L^-)^{1/2} \zeta , \zeta \rangle = \langle L^+ \xi , \xi \rangle <0$$
with $\xi:=(L^-)^{1/2} \zeta$. Let $P^-$ be the orthogonal projection from $L^2$ to $\text{Ker}(L^-)^\perp$.
Thanks to the Rayleigh principle, the operator $P^- L^+ P^-$ has a negative eigenvalue
$E_3 \in [-E^*,0)$: $L^+ \psi = E_3 \psi + c \phi_\mu$ for some $\psi \perp \phi_\mu$, $\psi \neq 0$, $c \in \mathbb{C}$.
If $c=0$, then $\psi$ is the ground state of $L^+$ thus $\psi>0$;
in particular the two positive functions $\psi$ and $\phi_\mu$ cannot be orthogonal in $L^2(0,1)$: contradiction.
Thus, $c \neq 0$ and $(L^+-E_3)^{-1} \phi_\mu = \frac{\psi}{c}$. In particular, we have
$$g(E_3)=\langle (L^+-E_3)^{-1} \phi_\mu, \phi_\mu \rangle = \frac{1}{c} \langle \psi , \phi_\mu \rangle =0.$$
which is impossible in view of (\ref{absurd}). Therefore, the eigenvalues of $\mathcal{T}$ are $\geqslant 0$.
\\

\noindent \emph{Third step: $\text{Sp}(\mathcal{T}) \subset \mathbb{R}^+$ in the defocussing case.}
Let us assume that $\mathcal{T}$ has a negative eigenvalue $E<0$.
Let $\psi$, $\psi_1$, $c$ be as in the first step and $\xi:=L^- \psi_1$. Then
$$\begin{array}{ll}
0 & < \langle L^+ \xi , \xi \rangle  = \langle L_\mu^+ L^- \psi_1 , L^- \psi_1 \rangle \\
  & = \langle E(\psi_1+c\phi_\mu) , L^- \psi_1 \rangle \\
  & = E  \| (L^-)^{1/2} \psi_1 \|_{L^2}^2 \quad \text{ because } L^- \phi_\mu =0 \\
  & < 0 : \text{contradiction.}
\end{array}$$
Therefore, the eigenvalues of $\mathcal{T}$ are $\geqslant 0$.
\\

\noindent \emph{Fourth step: Conclusion.} Thanks to (\ref{Lop2}) and the second and third steps, 
the eigenvalues of $\Lop^2$ are non positive real numbers.
Thus, the eigenvalues of $\Lop$ are purely imaginary. $\hfill \Box$

\subsection{Statement \emph{(iv)}}

Note that $0$ is an eigenvalue of $\Mop_{\mp \pi^2}$ with multiplicity $2$:
$$\Mop_{\mp \pi^2} \left( \begin{array}{c}
\sin(\pi x) \\ 0
\end{array}\right)=0, \qquad
\Mop_{\mp \pi^2} \left( \begin{array}{c}
0 \\ \sin(\pi x) 
\end{array}\right)=0$$
and the non zero eigenvalues of $\Mop_{-\pi^2}$ are $\{ \pm (n^2-1)\pi^2 ; n \geqslant 2 \}$:
$$\Mop_{\mp \pi^2} \left( \begin{array}{c}
\sin(n \pi x) \\ 0
\end{array}\right)=
(n^2-1)\pi^2 \left( \begin{array}{c}
\sin(n \pi x) \\ 0
\end{array}\right),$$ 
$$\Mop_{\mp \pi^2} \left( \begin{array}{c}
0 \\ \sin(n \pi x) 
\end{array}\right)=
-(n^2-1)\pi^2\left( \begin{array}{c}
0 \\ \sin(n \pi x) 
\end{array}\right).$$
For $\mu_0 \in [\mp \pi^2,+\infty)$, $\Mop_\mu$ converges to $\Mop_{\mu_0}$ when $\mu \rightarrow \mu_0$
in the sense of the generalized convergence of closed operators
(i.e. convergence of the graph, see \cite[Chapter IV, paragraph 2, page 197]{K66}).
Thus $\mu \mapsto \beta_{n,\mu}$ is continuous for every $n \in \mathbb{N}$
(see \cite[Chapter IV, paragraph 3.5]{K66}). $\hfill \Box$

\subsection{Statement \emph{(v)}}

In this section, we omit $\mu$ in subscript to simplify the notations.
Let $n \in \mathbb{N}^*$. The map
$$\begin{array}{|ccc}
\text{Ker}(\Lop - i \beta_{n} \text{Id}) & \rightarrow & \mathbb{C}^2 \\
               \Phi                      &  \mapsto    &  \Phi'(0)
\end{array}$$
is injective thanks to the uniqueness in Cauchy-Lipschitz theorem. Thus
$\text{dim}[\text{Ker}(\Lop - i \beta_{n} \text{Id})] \leqslant 2$.
\\

Now, we prove by contradiction that no non zero eigenvalue possesses a generalized eigenvector.
This proof follows the one of \cite{RSS05} for NLS on the whole space.
We use the operator $\mathcal{T}$ introduced in (\ref{Lop2}).
We assume that $\Lop_\mu$ has a generalized eigenvector associated to a non zero eigenvalue.
\\

\noindent \emph{First step: $\mathcal{T}$ has a generalized eigenvector associated to a non zero eigenvalue.}
Let $\psi, \rho \in D(\Lop_\mu)-\{0\}$ and $E \neq 0$ be such that
$$\Lop_\mu \rho = E \rho, \qquad \Lop_\mu \psi=E\psi + \rho.$$
Then,
$$(\Lop_\mu^2-E^2)\rho =0, \qquad (\Lop_\mu^2-E^2) \psi=2E\psi,$$
thus $\Lop_\mu^2$ has a generalized eigenvector associated to the eigenvalue $2E$, and so has $\mathcal{T}$ (see \ref{Lop2}).
\\

\noindent \emph{Second step: $(L^-)^{1/2} L^+ (L^-)^{1/2}$ has a generalized eigenvector.}
Let $\psi, \rho \in D(\mathcal{T})-\{0\}$, $E \neq 0$ be such that
$$\mathcal{T}\psi=E\psi, \qquad \mathcal{T}\rho=E\rho+c\psi.$$
Then $\psi_1$ and $\rho_1$ are not collinear to $\phi_\mu$ (otherwise $E$ would be zero).
Let $\psi_1, \rho_1$ be the projections of $\psi, \rho$ orthogonally to $\phi_\mu$. Then 
$(L^-)^{1/2} \psi_1 \neq 0$, $(L^-)^{1/2} \rho_1 \neq 0$ (because of (\ref{Ker_L})) and
$$[(L^-)^{1/2} L^+ (L^-)^{1/2}-E] (L^-)^{1/2} \psi_1 =0
\qquad
[(L^-)^{1/2} L^+ (L^-)^{1/2}-E] (L^-)^{1/2} \rho_1 = c (L^-)^{1/2} \psi_1.$$
thus $(L^-)^{1/2} L^+ (L^-)^{1/2}$ has a generalized eigenvector.
\\

\noindent \emph{Third step: The operator $B:=(L^-)^{1/2} L^+ (L^-)^{1/2}$ with domain
$D(B):=H^4_{(0)}(0,1)$ is self adjoint, which gives the contradiction.} The symmetry of $B$ is obvious.
Let us prove that $D(B^*)=D(B)$. Let $g,h \in L^2(0,1)$ be such that 
$$\langle (L^-)^{1/2} L^+ (L^-)^{1/2} f , g \rangle = \langle f , h \rangle, \quad \forall f \in H^4_{(0)}(0,1).$$
Our goal is to prove that $g \in H^4_{(0)}(0,1)$.
Taking $f\in \text{Ker}[(L^-)^{1/2}]$ shows that $h \in \text{Ker}[(L^-)^{1/2}]^{\perp}$.
By Fredholm alternative applied to the self adjoint operator $(L^-)^{1/2}$,
there exists $h_1 \in D[(L^-)^{1/2}]=H^1_0(0,1)$ such that $h=(L^-)^{1/2}h_1$. Then
$$\langle (L^-)^{1/2} L^+ (L^-)^{1/2} f , g \rangle = \langle f , (L^-)^{1/2} [ h_1 + c \phi_\mu] \rangle, \quad \forall f \in H^4_{(0)}(0,1), \quad \forall c \in \mathbb{C}.$$
By self-adjointness of $(L^-)^{1/2}$ and (\ref{Ker_L}), this gives
$$\langle (L^-)^{1/2} L^+ f_1 , g \rangle = \langle f_1 , h_1 + c \phi_\mu \rangle, \quad \forall f_1 \in H^3_{(0)}(0,1) \text{ with } f_1 \perp \phi_\mu, \quad \forall c \in \mathbb{C}.$$
The restriction '$f_1 \perp \phi_\mu$' may be removed by choosing
$$c:=\frac{1}{\|\phi_\mu\|_{L^2}^2} \Big( \langle(L^-)^{1/2} L^+ \phi_\mu,g\rangle - \langle\phi_\mu,h_1\rangle  \Big).$$
Then,
\begin{equation} \label{f1}
\langle (L^-)^{1/2} L^+ f_1 , g \rangle = \langle f_1 , h_1 + c \phi_\mu \rangle, \quad \forall f_1 \in H^3_{(0)}(0,1).
\end{equation}
Thanks to (\ref{Ker_L}), the operator $L^+:H^3_{(0)}(0,1) \rightarrow H^1_0(0,1)$ is bijective and selfadjoint thus
$$\langle (L^-)^{1/2} f_2 , g \rangle = \langle f_2 , (L^+)^{-1} [h_1 + c \phi_\mu] \rangle, \quad \forall f_2 \in H^1_0(0,1).$$
By selfadjointness of $(L^-)^{1/2}$, this proves that
$$(L^-)^{1/2}g= (L^+)^{-1} [h_1 + c \phi_\mu]$$
belongs to $H^3_{(0)}(0,1)$ (because $h_1 \in H^1_0(0,1)$), thus $g \in H^4_{(0)}(0,1)$. $\hfill \Box$

\subsection{Statements \emph{(vi)} and \emph{(vii)}}

One easily checks that (\ref{jordan}) holds.
\\

\noindent \emph{First step: $\text{Ker}(\Lop_\mu)=\mathbb{C} \Phi_0^+$, $\forall \mu \in (\mp \pi^2,+\infty)$.}
Let $(u,v) \in D(\Lop_\mu)$ be such that $L_\mu^-v=L_\mu^+u=0$.
From (\ref{Ker_L}), we deduce that $u=0$ and $v=c\phi_\mu$ for some $c \in \mathbb{R}$.
\\

\noindent \emph{Second step: $\Lop_\mu$ does not have a third (linearly independent) generalized eigenvector, for every $\mu \in (\mp \pi^2,+\infty)$.}
We assume that there exists $(u,v) \in D(\Lop_\mu)$ such that
$L_\mu^- v = \partial_\mu \phi_\mu$ and $L_\mu^+ u =0$. Then, thanks to (\ref{convexite}) and the selfadjointness of $L_\mu^-$, we get
$$0 < \langle \partial_\mu \phi_\mu , \phi_\mu \rangle = \langle L_\mu^- v , \phi_\mu \rangle 
= \langle v , L_\mu^- \phi_\mu \rangle = 0,  $$
which is impossible. 
\\

This proves that $(\Phi_0^-,\Phi_0^+)$ form a basis of the generalized null space for $\Lop_\mu$.
The case of $\Lop_\mu^*$ may be treated similarly. Moreover, we have
$$ \sigma i \beta_{m,\mu} \langle \Zn_m^\sigma ,\Psi _n^\tau \rangle =
\langle \Lop_\mu \Zn_m^\sigma , \psi_n^\tau \rangle =
\langle \Zn_m^\sigma , \Lop_\mu^* \psi_n^\tau \rangle =
\tau i \beta_{n,\mu} \langle \Zn_m^\sigma ,\Psi _n^\tau \rangle.$$
This proves (\ref{Kronecker}) when all the positive eigenvalues of $\Lop_\mu$ are simple.

\section{Analyticity of eigenvalues: proof}

The proof of Proposition \ref{Prop:vap_analytic} 
relies on the fact that the dimension of the eigenspaces of $\Mop_\mu$ is at most two, and the following elementary result.

\begin{prop} \label{Prop:matrice22}
Let $I \subset \mathbb{R}$ be an interval and
$B:I \rightarrow \mathcal{M}_2(\mathbb{R})$ be an analytic function.
Assume that the eigenvalues of $B(\mu)$ are real for every $\mu \in I$. Then, there exists analytic functions
$\lambda_1,\lambda_2:I \rightarrow \mathbb{R}$ such that
$\text{Sp}[B(\mu)]=\{ \lambda_1(\mu),\lambda_2(\mu)\}$ for every $\mu \in I$.
\end{prop}

\noindent \textbf{Proof of Proposition \ref{Prop:matrice22}:} The eigenvalues of $B(\mu)$ are
\begin{equation} \label{VA}
\lambda_{\pm}(\mu):=\frac{1}{2} \Big( \text{Tr}[B(\mu)] \pm \sqrt{\Delta(\mu)} \Big)
\text{ where }
\Delta(\mu):=\text{Tr}[B(\mu)]^2 - 4 \text{Det}[B(\mu)].
\end{equation}
Let $\mu_0 \in I$. If $\Delta(\mu_0)>0$, then the previous formula defines 2 analytic functions on a neighborhood of $\mu_0$.
Let us assume that $\Delta(\mu_0)=0$. Notice that $\Delta(\mu)\geqslant 0$, $\forall \mu \in I$ because $A(\mu)$ has real eigenvalues.
Expanding $\Delta(\mu)$ in power series of $(\mu-\mu_0)$, we find $k \in \mathbb{N}^*$ and an function
$\widetilde{\Delta}(\mu)$, analytic in a neighborhood of $\mu_0$ and satisfying $\widetilde{\Delta}(\mu_0)>0$ such that
$\Delta(\mu)=(\mu-\mu_0)^{2k} \widetilde{\Delta}(\mu)$ on a neighborhood of $\mu_0$. Then we get the conclusion with the formula
$$\lambda_1(\mu):= \frac{1}{2} \Big( \text{Tr}[B(\mu)] - (\mu-\mu_0)^k \sqrt{\widetilde{\Delta}(\mu)} \Big), \quad
\lambda_2(\mu):= \frac{1}{2} \Big( \text{Tr}[B(\mu)] + (\mu-\mu_0)^k \sqrt{\widetilde{\Delta}(\mu)} \Big). \hfill \Box $$

\begin{prop} \label{Prop:analytic_everywhere}
Let $\mu_0 \in (\mp \pi^2,\infty)$ and $n \in \mathbb{N}^*$.
There exists an analytic function $\varphi$, defined on an open neighborhood $I$ of $\mu_0$ such that
$\varphi(\mu_0)=\beta_{n,\mu_0}$ and $\varphi(\mu) \in \text{Sp}(\Mop_\mu)$, $\forall \mu \in I$.
\end{prop}

\noindent \textbf{Proof of Proposition \ref{Prop:analytic_everywhere}:} Let $\mu_0 \in (\mp \pi^2,\infty)$.
\\

\noindent \emph{First step: Reduction to a finite dimensional space.}
\\
This step follows exactly \cite[Chap VII, Paragraph 1.3, proof of Theorem 1.7, page 368]{K66}.
Let $\mathcal{C}$ be a closed curve in the complex plane 
that separates $\text{Sp}(\Mop_{\mu_0})$ into two parts:
a finite one $\Sigma'(\mu_0)$, with cardinal $N \in \mathbb{N}^*$ and an infinite one $\Sigma''(\mu_0)$.
Since $\Mop_\mu$ converges to $\Mop_{\mu_0}$ when $\mu \rightarrow \mu_0$, in the generalized sense 
(convergence of graphs of closed operators), then, for sufficiently small $|\mu-\mu_0|$,
$\text{Sp}(\Mop_{\mu})$ is likewise separated by $\mathcal{C}$ into a finite part
$\Sigma'(\mu)$, with cardinal $N$, and an infinite part $\Sigma''(\mu)$,
associated to the decomposition $L^2((0,1),\mathbb{C}^2)=E'(\mu) \oplus E''(\mu)$.
The projection on $E'(\mu)$ along $E''(\mu)$ is given by
$$P(\mu)=\frac{1}{2\pi i} \int_{\mathcal{C}} (\Mop_\mu-z \text{Id})^{-1} dz.$$
It is a bounded-holomorphic operator near $\mu=\mu_0$.
\\

Let us construct a transformation $U(\mu)$ such that
\begin{enumerate}
\item $U(\mu)$ and $U(\mu)^{-1}$ are bounded-holomorphic on $L^2((0,1),\mathbb{C}^2)$,
\item $U(\mu) P(\mu_0) U(\mu)^{-1}=P(\mu)$ for every $\mu$ near $\mu_0$.
\end{enumerate}
We define $U(\mu)$ and $V(\mu)$ as the operators on $L^2((0,1),\mathbb{C}^2)$,
solutions of the linear ordinary differential equations
$$\begin{array}{cc}
\left\lbrace \begin{array}{l}
U'(\mu)=Q(\mu) U(\mu),\\
U(\mu_0)=\text{Id},
\end{array}\right.
&
\left\lbrace \begin{array}{l}
V'(\mu)= - V(\mu)Q(\mu),\\
V(\mu_0)=\text{Id},
\end{array}\right.
\end{array}$$
where $Q(\mu):=P'(\mu)P(\mu)-P(\mu)P'(\mu)$. Then, $U(\mu)$ and $V(\mu)$ are bounded-holomorphic and
$$(VU)'=V'U+VU'=-VQU+VQU=0$$
thus $V(\mu)U(\mu) \equiv \text{Id}$. This proves the announced properties on $U(\mu)$.
\\

Note that $$\hat{\Mop_\mu}:=U(\mu)^{-1} \Mop_\mu U(\mu)$$
commutes with $P(\mu_0)$. Indeed, $\Mop_\mu$ commutes with $P(\mu)$ thus the property (ii) above proves 
$$\begin{array}{ll}
\hat{\Mop_\mu} P(\mu_0)
& =U(\mu)^{-1} \Mop_\mu P(\mu)U(\mu) = U(\mu)^{-1} P(\mu) \Mop_\mu U(\mu)\\
& =P(\mu_0)U(\mu)^{-1}\Mop_\mu U(\mu) = P(\mu_0)\hat{\Mop_\mu}.
\end{array}$$
Thus, the N-dimensional space $E'(\mu_0)=\text{Range}[P(\mu_0)]$ is stable by $\hat{\Mop_\mu}$ and
\begin{equation} \label{reduc}
\text{Sp}[\hat{\Mop_\mu}|_{E'(\mu_0)}]=\Sigma'(\mu).
\end{equation}
\\

\noindent \emph{Second step: Analyticity of eigenvalues.}
\\
Let $n \in \mathbb{N}^*$. We apply the first step with a positively oriented circle $\mathcal{C}$ with center $\beta_{n,\mu_0}$
and radius $\epsilon>0$ small enough so that $\mathcal{C}$ contains no other eigenvalue of $\Mop_{\mu_0}$.
If $\beta_{n,\mu_0}$ is simple, then the previous construction shows that $\mu \mapsto \beta_{n,\mu}$ is analytic near $\mu=\mu_0$.
Let us assume that $\beta_{n,\mu_0}$ is a multiple eigenvalue of $\Mop_{\mu_0}$.
Thanks to Proposition \ref{prop:basic_spect} \emph{(v)}, $E'(\mu_0):=\text{Ker}[\Mop_{\mu_0}-\beta_{n,\mu_0}\text{Id}]$ has dimension $2$.
Let $(e_1,e_2)$ be a basis of $E'(\mu_0)$. One may assume that $e_1$ and $e_2$ are real-valued functions, otherwise consider $(e_j+\overline{e_j})/2$ and
$(e_j-\overline{e_j})/(2i)$. Let $B(\mu)$ be the $2*2$-matrix of the operator $\hat{\Mop_\mu}|_{E'(\mu_0)}$ in the basis $(e_1,e_2)$.
Then $B(\mu)$ is analytic and has only real valued eigenvalues, thanks to (\ref{reduc}) and Proposition \ref{prop:basic_spect} \emph{(ii)}.
Let us prove that $B(\mu)$ has real valued coefficients, which allows to conclude thanks to Proposition \ref{Prop:matrice22}.
\\

\emph{Step 2.1: We prove that $P(\mu)$ is real valued, i.e. $P(\mu)f \in L^2((0,1),\mathbb{R}^2)$, $\forall f \in L^2((0,1),\mathbb{R}^2)$.}
Indeed,
\begin{eqnarray*}
P(\mu)f & =& \frac{1}{2\pi i} \int_{\mathcal{C}} (\Mop_\mu-z \text{Id})^{-1} f dz =
\frac{1}{2\pi} \int_0^{2\pi}  \Big( \Mop_\mu- (\beta_{n,\mu_0}+\epsilon e^{i\theta}) \text{Id} \Big)^{-1} f \epsilon e^{i\theta} d\theta \\
        & = &\overline{\frac{1}{2\pi} \int_0^{2\pi}  \Big( \Mop_\mu- (\beta_{n,\mu_0}+\epsilon e^{-i\theta}) \text{Id} \Big)^{-1} f \epsilon e^{-i\theta} d\theta} \\
        & = &\overline{\frac{1}{2\pi i} \int_{\mathcal{C}} (\Mop_\mu-z \text{Id})^{-1} f dz}  = \overline{P(\mu)f}.
\end{eqnarray*}

\emph{Step 2.2: We prove that $U(\mu)$ and $U(\mu)^{-1}$ are real valued.} Indeed, if $f \in L^2((0,1),\mathbb{R}^2)$, then
$g(\mu):=U(\mu)f$ and is the solution of the ordinary differential equation
$$\left\lbrace \begin{array}{l}
g'(\mu)=Q(\mu) g(\mu),\\
g(\mu_0)=f,
\end{array}\right.$$
thus it is real valued, thanks to Step 2.1 and the uniqueness in Cauchy-Lipschitz theorem.
\\

\emph{Step 2.3: We prove that $B(\mu)$ have real valued coefficients.} 
Thanks to \emph{Step 2.2}, we have
$$B(\mu)e_j = \hat{\Mop_\mu} e_j = U(\mu)^{-1} \Mop_\mu U(\mu) e_j \in L^2((0,1),\mathbb{R}^2), \quad \forall j=1,2.$$
thus its coefficients on the (real-valued) basis $(e_1,e_2)$ are real. $\hfill \Box$
\\
\\

\noindent \textbf{Proof of Proposition \ref{Prop:vap_analytic}:}
By \cite[Chapter VII, paragraph 3, Theorem 1.8]{K66}, the eigenvalues of $\Mop_\mu$ are branches of one or several
analytic functions, which have only algebraic singularities, and which are everywhere continuous. An exceptional point $\mu_0$ is
\begin{itemize}
\item either a branch point (see \cite[Chap II, Paragraph 1.2]{K66} for a definition),
\item or a regular point where different eigenvalues coincide (crossing).
\end{itemize}
Moreover, when we consider a finite number of eigenvalues,
there are only a finite number of exceptional points $\mu_0$ in each compact set of $(\mp \pi^2,\infty)$.
Proposition \ref{Prop:analytic_everywhere} shows that there are no branch point and that eigenvalues can be followed analytically through crossings.
\\

Let $n \in \mathbb{N}^*$. There exists $\delta>\mp \pi^2$ such the map $\mu \mapsto \beta_{n,\mu}$ is continuous on $[\mp \pi^2,\delta)$,
and $\beta_{n,\mu}$ is a simple eigenvalue of $\Mop_\mu$ for every $[\mp \pi^2,\delta)$.
Then, $\mu \mapsto \beta_{n,\mu}$ is analytic on $(\mp \pi^2,\delta)$ thanks to Proposition \ref{Prop:analytic_everywhere}.
Let $\mu^*$ be the sup of the $\mu_\sharp \geqslant \delta$ such that 
$\mu \in (\mp \pi^2,\delta) \mapsto \beta_{n,\mu}$ may be extended in an analytic function $\varphi:(\mp \pi^2,\mu_\sharp) \rightarrow \mathbb{R}$,
which is everywhere an eigenvalue of $\Mop_\mu$. We prove by contradiction that $\mu^*=\infty$.

We assume that $\mu^*<+\infty$. Then at most a finite number of crossings may happen on $(\mp \pi^2,\mu^*)$:
there exists a finite number $N \in \mathbb{N}$ of points $\mu_1,...,\mu_N \in (\delta,\mu^*)$ such that
$\varphi(\mu)$ coincide with different eigenvalues
$\beta_{n_{k-1},\mu}$ when $\mu<\mu_k$ and $\beta_{n_k,\mu}$ when $\mu>\mu_k$,
with $n_k=n_{k-1}\pm 1$, for $k=1,...,N$.
In particular, for $\mu \in (\mu_N,\mu^*)$, we have $\varphi(\mu)=\beta_{n_N,\mu}$.
Thanks to Proposition \ref{Prop:analytic_everywhere}, 
$\varphi(\mu)$ may be extended into an analytic function on a larger interval than $(\mp \pi^2,\mu^*)$,
that is everywhere an eigenvalue of $\Mop_\mu$, which is impossible. 
Therefore $\mu^*=\infty$ and Proposition \ref{Prop:vap_analytic} is proved. $\hfill \Box$

\section{Moment problem}

The following proposition is crucial in the controllability of the linearized system.
It is a consequence of the Ingham inequality proved by Haraux in \cite{Haraux}
and may be proved exactly as \cite[Corollary 2 in Appendix B]{KB-CL}.

\begin{prop} \label{prop:moment_pb}
Let $T>0$, $N \in \mathbb{N}$ and $(\omega_{k})_{k \in \mathbb{N}}$ be an increasing sequence of $(0,+\infty)$ such that 
$\omega_{k+1} - \omega_{k} \rightarrow + \infty$ when $k \rightarrow + \infty$.
Then there exists a continuous linear map
$$\begin{array}{|cccc}
L: & \mathbb{R}^N \times l^2(\mathbb{N}^*,\mathbb{C}) & \rightarrow & L^2((0,T),\mathbb{R}) \\
   &        (\tilde{d},d)                           & \mapsto     & L(\tilde{d},d)
\end{array}$$
such that, for every $\tilde{d}=(\tilde{d}_1,...,\tilde{d}_N) \in \mathbb{R}^N$ and
$d=(d_k)_{k \in \mathbb{N}} \in l^2(\mathbb{N}^*,\mathbb{C})$,
the function $v:=L(\tilde{d},d)$ solves
$$\begin{array}{l}
\int_{0}^{T} v(t) e^{i \omega_{k} t} dt = d_{k} , \forall k \in \mathbb{N}^*, \\
\int_0^T t^k v(t) dt = \tilde{d}_{k+1}, \forall k=0,...,N-1.
\end{array}$$
\end{prop}


\bibliographystyle{plain}

\end{document}